\documentclass[preprint,10pt]{elsarticle}

\usepackage{amsmath,amssymb,amsthm,mathtools,graphicx,subfigure,version}
\usepackage[colorlinks,
  linkcolor = black,
  anchorcolor = black,
  citecolor = black]{hyperref}
\usepackage[nameinlink,capitalize]{cleveref}
\usepackage{booktabs,version}
\usepackage[bottom,stable]{footmisc}
\usepackage{graphicx}
\usepackage{tikz}
\usepackage{algorithm}
\usepackage{algpseudocode}
\numberwithin{algorithm}{section}
\algrenewcommand\algorithmicrequire{\textbf{Input:}}
\algrenewcommand\algorithmicensure{\textbf{Output:}}

\def\re{{\rm e}}
\def\rd{{\rm d}}
\def\ri{{\rm i}}
\newcommand{\bs}[1]{\boldsymbol{#1}}
\graphicspath{{./figs/}}

\newtheorem{remark}{Remark}[section]

\numberwithin{equation}{section}
\numberwithin{table}{section}
\numberwithin{figure}{section}

\journal{}

\begin{document}
\begin{frontmatter}

\title{Lightning Algorithm for Inhomogeneous Modified Helmholtz Equations on Cornered Domains and Their Applications\tnoteref{funding}}
\tnotetext[funding]{This work is supported by the National Natural Science Foundation of China under Grants 12601710, W2431008, 12531015 and 12371409.}

\author[eitech]{Desong Kong}
\ead{kongdesong2022@outlook.com}
\author[eitech]{Jie Shen\corref{cor1}}
\ead{jshen@eitech.edu.cn}
\cortext[cor1]{Corresponding author.}
\address[eitech]{School of Mathematical Sciences, Eastern Institute of Technology, Ningbo 315200, China}

\begin{abstract}
We develop in this paper an efficient solver for the inhomogeneous modified Helmholtz equations in two-dimensional complex  domains with particular attention to resolving corner singularities.
We represent the solution as a summation of a particular solution and the solution of the modified Helmholtz equation.
Built upon suitable domain embedding methods, we compute the particular solution using the Fourier spectral method and/or Legendre--Fourier spectral method.
To treat the corner singularities, we develop a lightning algorithm by employing an expansion of the modified Bessel functions of the second kind, coupled with tapered exponentially clustering poles, to approximate the corner singularities, and using an expansion of the modified Bessel functions of the first kind to approximate the smooth part of the solution.
Finally, we determine the unknown expansion coefficients by solving a least-squares problem to satisfy the boundary conditions.
Numerical experiments demonstrate that root-exponential convergence can be achieved  for the modified Helmholtz equation  with corner singularities.
To further demonstrate the applicability of the proposed algorithm, it has been successfully applied to  solve some typical linear and nonlinear  parabolic equations. 
We have also made a first attempt in extending the lightning algorithm to a three-dimensional modified Helmholtz equation.

\end{abstract}


\end{frontmatter}

\section{Introduction}\label{sec:intro}
We develop in this paper a fast and accurate numerical algorithm for solving  the  problem
\begin{align}\label{eq:prob1}
    -\Delta u+\kappa^2 u=f \quad \text{in }\Omega;
    \quad u=h \quad \text{on}\;\; \partial\Omega,
\end{align}
where $\kappa \in \mathbb R$, $\Omega\subset\mathbb R^d$ ($d=2,3$) is a bounded domain whose boundary  has corners or singular points, and $f$ and $h$ are given functions. Problems of this type arise repeatedly after implicit or implicit-explicit time discretizations of parabolic equations such as
\begin{equation*}
    \partial_t u = \nu\Delta u + \mathcal N(u,t),
    \quad \text{in}\;\; \Omega,\quad t\in(0,T],
\end{equation*}
with suitable initial and boundary conditions. An efficient solver for \eqref{eq:prob1} is therefore desiable  as a reusable kernel in time-dependent computations or as a stand-alone elliptic solver.

Solutions of elliptic problems on cornered domains generally have reduced regularity near the corners. Uniform-grid finite differences, finite elements, and standard spectral methods may therefore require substantial local refinement or geometric decomposition; see, for example, \cite{ciarlet2002finite,brenner2008mathematical,langtangen2017solving,Trefethen2000Book,Canuto2006Book,Shen2011Book,mortensen_joss}. Lightning methods address this difficulty by placing singular basis functions at preassigned exterior poles that cluster exponentially near the corners. The original method was developed for Laplace problems \cite{gopal2019solving}, and related solvers have since been constructed for Helmholtz and Stokes equations \cite{gopal2019new,brubeck2022lightning}. Other developments include lightning potential-flow solvers \cite{baddoo2020lightning}, log-lightning computation of logarithmic capacity and Green's functions \cite{baddoo2021loglightning}, and the AAAtrig method for periodic rational approximation \cite{baddoo2021aaatrig}.

A related modified-Bessel lightning representation has recently been used for heat-equation computations \cite{croix2025lightn}, while our emphases are on inhomogeneous modified Helmholtz problems \eqref{eq:prob1} and its applications in solving nonlinear parabolic equations.

For $f\not\equiv0$, write $u=u_p+v$, where a particular solution satisfies
\begin{equation}\label{eq:sol-part}
    -\Delta u_p+\kappa^2 u_p=f \quad \text{in}\;\; \Omega,
\end{equation}
without a prescribed boundary condition. The correction $v$ then solves
\begin{align}\label{eq:prob10}
    -\Delta v+\kappa^2 v=0 \quad \text{in}\;\; \Omega,
    \quad
    v=g:=h-u_p|_{\partial\Omega} \quad \text{on}\;\; \partial\Omega.
\end{align}
For $\kappa \neq 0$, we approximate $v$ by modified-Bessel multipoles centered at exterior poles and by a regular expansion centered inside the domain. The case $\kappa=0$ is treated by the original Laplace lightning algorithm. To compute $u_p$, we extend $f$ to a regular enclosing domain and solve the extended equation by a fast spectral method with either a Fourier-Legendre spectral method on a disk or  a Legendre approach on a rectangle.

The principal contributions of this paper are as follows.
\begin{itemize}
    \item[(i)] To deal with the homogeneous modified Helmholtz equation (i.e., \eqref{eq:prob10} with $f=0$) in 2D, we use the modified Bessel function of the second kind, combined with the tapered exponentially clustered poles \eqref{eq:pole-trap},  to capture the singular part of the solution, and  
     we employ the modified Bessel functions of the first kind to approximate the regular part of the solution, leading to root-exponential convergence for the homogeneous modified Helmholtz equation. 
     
    \item[(ii)] To deal with the non-homogeneous modified Helmholtz equation, we combine the homogeneous correction with a  particular solution by an embedded-domain method.
    
     \item[(iii)] To deal with nonlinear parabolic equations,  we combine  the modified Helmholtz  solver with a linearly implicit BDF type discretization.
    
    \item[(iv)] We extend the proposed method to three-dimensional problems with point singularities by introducing an ansatz that combines spherical harmonic expansions with modified Bessel functions.

\end{itemize}

The remainder of the paper is organized as follows. In \cref{sec:lig-2d}, we introduce the two-dimensional homogeneous lightning algorithm to the homogeneous modified Helmholtz equation. In \cref{sec:inhom-poi}, we introduce the embedded-domain method to compute the particular solution. Then in \cref{sec:app}, we present the truncation expansion based on the RBF frame and apply the combination algorithm to parabolic equations. In \cref{sec:opt-gam}, we discuss parameter selection and a one-dimensional prototype approximation. In \cref{sec:helm-3d}, we describe the extension to three-dimensional problems with point singularities, followed by concluding remarks in \cref{sec:conclud}.

\section{The lightning algorithm for homogeneous modified Helmholtz equations}\label{sec:lig-2d}
We  develop in this section a lightning algorithm for homogeneous modified Helmholtz equations in a domain bounded by a piecewise smooth Jordan curve with corners. For the sake of brevity, the presentation and all numerical tests in this paper are concerned with Dirichlet boundary conditions. Neumann and Robin conditions can in principle be treated by applying the corresponding boundary operator to every basis function, but their compatibility, uniqueness, and error estimates differ from the Dirichlet case and are not analyzed here.

The term `lightning algorithm' refers to the use of preassigned exterior poles that cluster exponentially near geometric singularities. As explained in \cite{trefethen2021exponential}, the terminology reflects the same concentration of fields that makes lightning preferentially strike sharp objects.

\subsection{A brief review  of the lightning  algorithm for Laplace  equations}
The lightning  algorithm, originally introduced by Gopal and Trefethen in \cite{gopal2019solving}, is an efficient and robust method for solving Laplace equations defined on a two-dimensional domain $\Omega$, which is bounded by a piecewise smooth Jordan curve with corners.
For a real-valued homogeneous Laplace solution, the basic approximation is
\begin{align}\label{eq:rn-app}
u(z)\approx \Re \{r_N(z)\}, \quad
r_N(z)=\sum_{j=1}^{N_1}\frac{a_j}{z-z_j} + \sum_{j=0}^{N_2}b_j(z-z_*)^j.
\end{align}
Here $z=x+\ri y$, the poles  $z_j\in\mathbb C\setminus\overline\Omega$, and $z_*\in\Omega$. The rational part is often called the {\it Newman} part and resolves corner singularities, while the polynomial {\it Runge} part approximates the regular component. The complex approximation has $N=N_1+N_2+1$ coefficients. 
Typically, the process of calculating a rational approximation with unknown poles $\{z_j\}_{j=1}^{N_1}$ involves nonlinear iterations \cite{hochman2013use}, which can be challenging and time-consuming. However, the scheme in \eqref{eq:rn-app} linearizes the rational approximation by fixing the poles in advance. 
In other words, since $r_N(z)$ in \eqref{eq:rn-app} satisfies the homogeneous Laplace equation, with ${z_j}$ and $z_*$ fixed, we  can enforce the boundary condition  on a discrete set of points on the boundary to  obtain a  least-squares problem to  determine $\{a_j\}$ and $\{b_j\}$.

\subsection{Lightning algorithm for two-dimensional homogeneous modified Helmholtz equations}\label{subsec:helm-2d}

Consider first the homogeneous problem~\eqref{eq:prob10} with $\kappa \neq 0$.  Let $z-z_*=r_{\! *}\, \re^{\ri\theta_*}$ for a fixed expansion center $z_*$.  
A local polar separation of variables gives
\begin{equation}\label{eq:four-series-2d}
  v(r_*,\theta_*)=\sum_{m=-\infty}^{\infty} \widehat v_m(r_*)\, \re^{\ri m\theta_*},
\end{equation}
and with $\eta=\kappa\,r_*$, each radial coefficient satisfies the modified Bessel equation
\begin{equation}\label{eq:mod-bessel}
  \eta^2\widehat w_m''(\eta)+\eta\widehat w_m'(\eta) - (\eta^2+m^2)\widehat w_m(\eta)=0,
  \quad m\in\mathbb Z,
\end{equation}
whose solutions are  linear combinations of  the modified Bessel functions of the first and second kinds, $I_m$ and $K_m$.  
Thus, %
the formal multipole expansion of $v$ is
\begin{equation}\label{eq:poi-formal}
  v(z) = \sum_{m=-\infty}^{\infty}
  \bigl(a_mK_{|m|}(\kappa\,r_*) + b_mI_{|m|}(\kappa\,r_*)\bigr) \, \re^{\ri m\theta_*}.
\end{equation}

We recall that, for an integer $m\geq1$ and $|s|\ll 1$,
\begin{equation}\label{eq:asymp-IK}
  I_m(s) \sim \frac{(s/2)^m}{m!},\quad
  K_m(s) \sim \frac{(m-1)!}{2}\Big(\frac{s}{2}\Big)^{-m},
  \quad K_0(s)\sim-\log s,
\end{equation}
while $I_{-m}=I_m$ and $K_{-m}=K_m$; see
\cite[(10.27.1), (10.27.3), and (10.30.1)--(10.30.3)]{Olver2010Book}.
Obviously, the $K_m$ terms are singular at their centers, whereas the $I_m$ terms are regular.

Motivated by the Newman-Runge decomposition in the Laplace lightning algorithm, we place the singular centers outside the physical domain and use a regular expansion about an interior point. 
To be more precise, we set
\begin{align*}
  r_n=|z-z_n|,\quad \theta_n=\arg(z-z_n),\quad
  r_*=|z-z_*|,\quad \theta_*=\arg(z-z_*),
\end{align*}
where $z_n\in\mathbb C\setminus\overline\Omega$ are exterior poles clustered near the corners and $z_*\in\Omega$.  
Hence, our computational ansatz can be written as
\begin{equation}\label{eq:poi-app}
\begin{split}
  R_N(z) &= \sum_{n=1}^{N_1}a_n K_1(\kappa\,r_n) \cos(\theta_n)
          + \sum_{n=1}^{N_1}b_n K_1(\kappa\,r_n) \sin(\theta_n) \\
          &\;\; + \sum_{m=0}^{N_2} c_mI_{m}(\kappa\,r_*) \cos(m\theta_*) 
          + \sum_{m=1}^{N_2} d_m I_{m}(\kappa\,r_*) \sin(m\theta_*).
\end{split}
\end{equation}
In this case, the number of unknowns is $N=2N_1+2N_2+1.$
Note that every basis function in~\eqref{eq:poi-app} satisfies $(-\Delta+\kappa^2)R_N=0$ in $\Omega$.  
The $K_1$ angular mode is the dipole-type field, while the $I_m$ block represents the regular component.  

\begin{remark}\label{rem:bessel}
It should be noted that the ansatz~\eqref{eq:poi-app} is not unique.  Indeed, one may include $K_0$ or higher-order multipoles in the exterior block, or replace part of the regular block by $K_m$ functions centered outside $\overline\Omega$.  
Recently, related modified-Bessel lightning representations have been used for heat-equation computations in \cite{croix2025lightn}.  
The present paper focuses on the bounded interior, inhomogeneous architecture described in \cref{sec:inhom-poi}-\cref{sec:app}.
\end{remark}

\subsubsection{Boundary least-squares formulation and residual control}
Let $v$ denote the real solution of \eqref{eq:prob10} and let
$e_N=v-R_N$. Since both $v$ and $R_N$ satisfy the homogeneous
modified Helmholtz equation in $\Omega$, we have
\[
  (-\Delta+\kappa^2)e_N=0 \qquad \text{in }\Omega.
\]
Assuming $e_N\in C^2(\Omega)\cap C(\overline{\Omega})$, the weak maximum principle for uniformly elliptic operators with nonnegative zeroth-order coefficient, applied to $e_N$ and $-e_N$, yields (cf., e.g., \cite[Sec.~3.1]{GilbargTrudinger2001} and \cite[Sec.~6.4.1]{evans2022partial})
\begin{equation}\label{eq:max-p}
\|e_N\|_{L^\infty(\Omega)} \leq \|e_N\|_{L^\infty(\partial\Omega)}.
\end{equation}
Thus, the continuous boundary residual provides an a posteriori $L^\infty$ bound for the error of the homogeneous correction.
Numerically, we can use a residual evaluated at finitely many boundary points to assess the accuracy of our proposed algorithm. More precisely,
let $\{\xi_j\}_{j=1}^M\subset\partial\Omega$ be an oversampled fitting grid with $M\geq N$, and set $g_j=g(\xi_j)$.  
With
\begin{align*}
  r_{j,n}=|\xi_j-z_n|,\quad \theta_{j,n}=\arg(\xi_j-z_n),\quad
  r_{j,*}=|\xi_j-z_*|,\quad \theta_{j,*}=\arg(\xi_j-z_*),
\end{align*}
the coefficient vector
\begin{align*}
  \bs c=(a_1,\ldots,a_{N_1},b_1,\ldots,b_{N_1},c_{0},\ldots,c_{N_2}, d_1, \ldots, d_{N_2})^\top \in\mathbb C^N
\end{align*}
can be determined from
\begin{equation}\label{eq:poi-least}
  \bs c=\operatorname*{argmin}_{\bs z\in\mathbb C^N}
  \|\bs A\bs z-\bs g\|_2,
\end{equation}
where $\bs g := (g_1, \ldots, g_M)^\top \in \mathbb C^M$, and the coefficient matrix $\bs A \in \mathbb C^{M \times N}$ has entries
\begin{equation*}\small
    \bs A_{j,n} = 
    \begin{dcases}
        K_{1}(\kappa\, r_{j,n})\; \cos(\theta_{j,n}), & 1 \le n \le N_1, \\
        K_{1}(\kappa\, r_{j,n})\; \sin(\theta_{j,n}), & N_1+1 \le n \le 2N_1, \\
        I_{n-2N_1-1}(\kappa\, r_{j,*})\; \cos((n-2N_1-1)\theta_{j,*}), & 2N_1+1 \le n \le N-N_2, \\ %
        I_{n-N+N_2}(\kappa\, r_{j,*})\; \sin((n-N+N_2)\theta_{j,*}), & N-N_2+1 \le n \le N,
    \end{dcases}
\end{equation*}

It should be noted that  the least-squares problem \eqref{eq:poi-least} is ill-conditioned due to  the exponentially clustered poles near the corners.  Therefore, a   suitable preconditioner should be used. 
For the rational approximation \eqref{eq:rn-app}, a technique based on the Arnoldi orthogonalization is applied to fix the exponentially ill-conditioned of the Vandermonde matrix \cite{brubeck2021vandermonde}.
While for the approximation \eqref{eq:poi-app}, we propose to use a technique based on the diagonal weighting technique (cf. \cite{Gopal2021}).
More precisely, we  use a diagonal preconditioner $\bs D$ which is obtained by normalizing the columns of $\bs A$, a diagonal preconditioner $\bs W$ with diagonal entries $\min_i|w_i-\xi_k|$ where $w_i$ is the location of the corner of $\Omega$, and solve the preconditioner least-squares problem 
\begin{align}\label{eq:column-scaling}
    \bs{\hat c} = \mathop{\rm arg min}_{\bs x \in \mathbb C^M} \| \bs W \bs A \bs D^{-1} \bs x - \bs W \bs g \|_2.
\end{align}
 The unknown coefficient vector $\bs c $ can be obtained by $\bs c = \bs D \bs{\hat c}.$

To assess whether the boundary sampling is sufficiently fine, we evaluate the residual on a second, independent validation grid $\bs \Xi \subset \partial \Omega$ that is not used in~\eqref{eq:poi-least}.
The reported validation residual is
\begin{equation*}%
  \mathcal E(\bs \Xi)
  =\max_{\bs \Xi} |R_N(\bs \Xi) - g(\bs \Xi)|.
\end{equation*}
In our following shown figures, we use a red circle to denote the value  $\mathcal E(\bs \Xi).$

\begin{figure}[!t]
    \centering
    \includegraphics[width=0.8\linewidth]{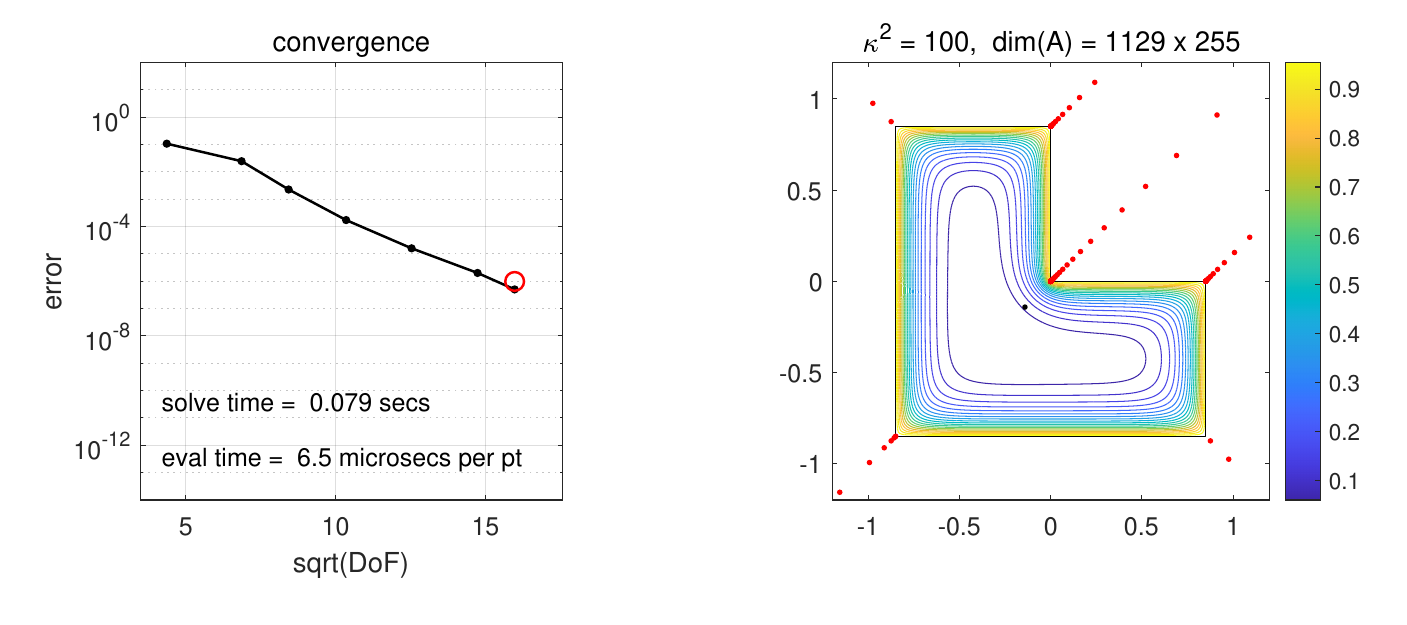}
    \vspace*{-10pt}
    \caption{Lightning solution of~\eqref{eq:prob10} on an L-shaped domain with
    $\kappa^2=100$ and $g=1$.  The convergence panel reports boundary fitting residuals
    and residuals on an independent validation grid.}
    \label{fig:poi-type-alp10}
\end{figure}

We first test~\eqref{eq:prob10} on an L-shaped domain with $\kappa^2=100$ and $g=1$.
The tapered exterior poles are defined in~\eqref{eq:pole-trap}, with the corner-dependent parameter \eqref{eq:optimal-sig}.  
\cref{fig:poi-type-alp10} shows the fitting and independent validation residuals together with the numerical solution.  
Their nearly linear decay against $\sqrt N$ is consistent with root-exponential convergence.  
The reported solve time is less than 0.1\,s in MATLAB on the test computer\footnote{All numerical experiments were performed in MATLAB R2023a on a desktop computer equipped with an Intel Core i5-13600KF CPU and 16 GB of RAM, running 64-bit Windows 11.}.
The figure title in the right panel means the value of $\kappa^2$ and the dimension of $\bs A$ whenever the algorithm achieves the preassigned tolerance $10^{-6}.$

\begin{figure}[!t]
    \centering
    \includegraphics[width=0.8\linewidth]{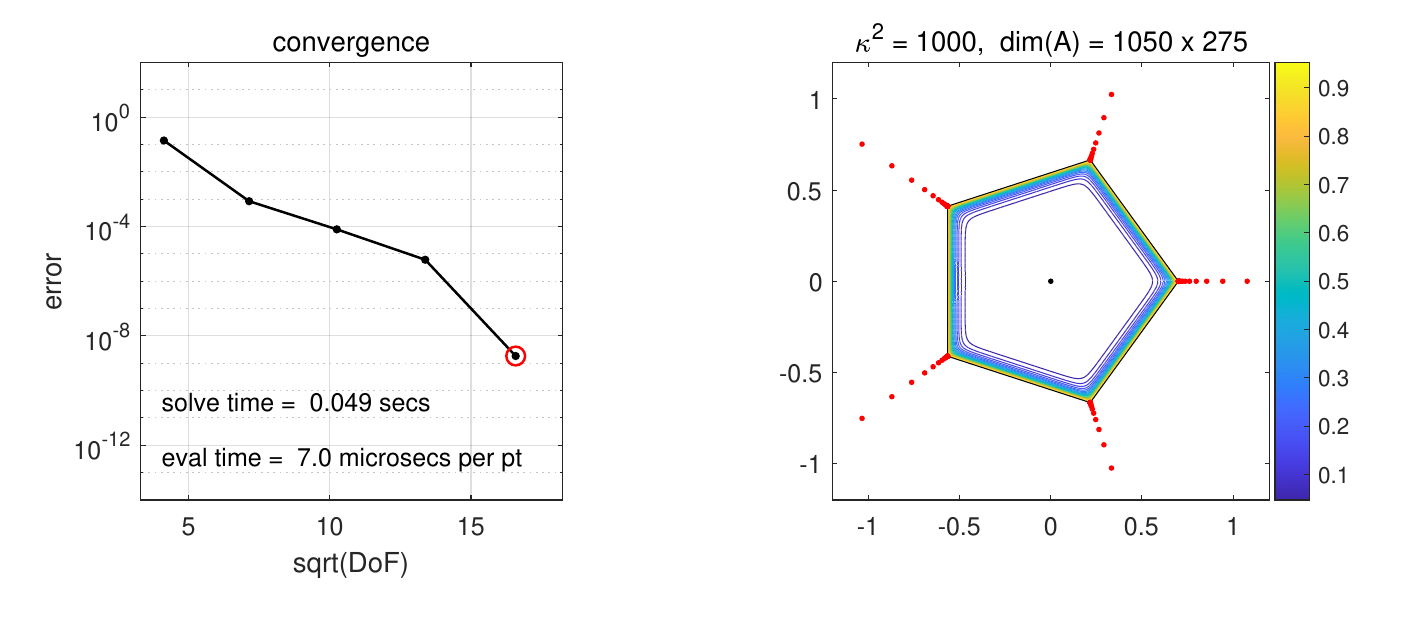}
    \includegraphics[width=0.8\linewidth]{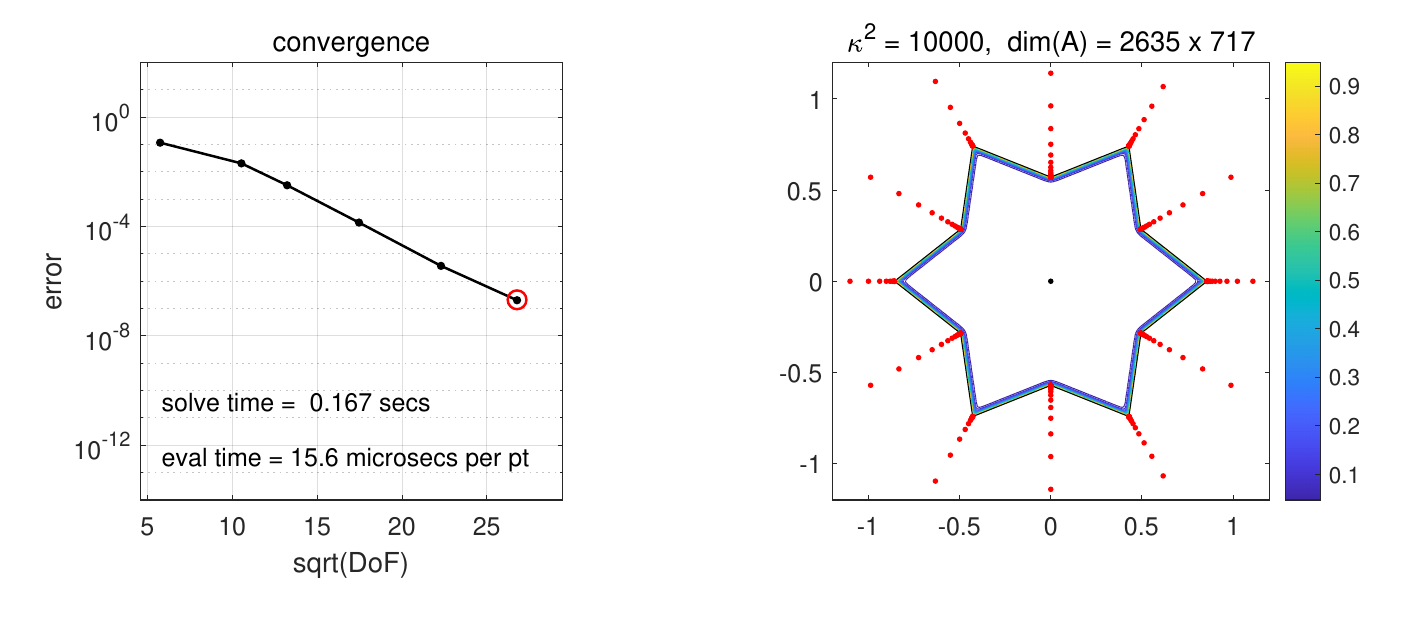}
    \includegraphics[width=0.8\linewidth]{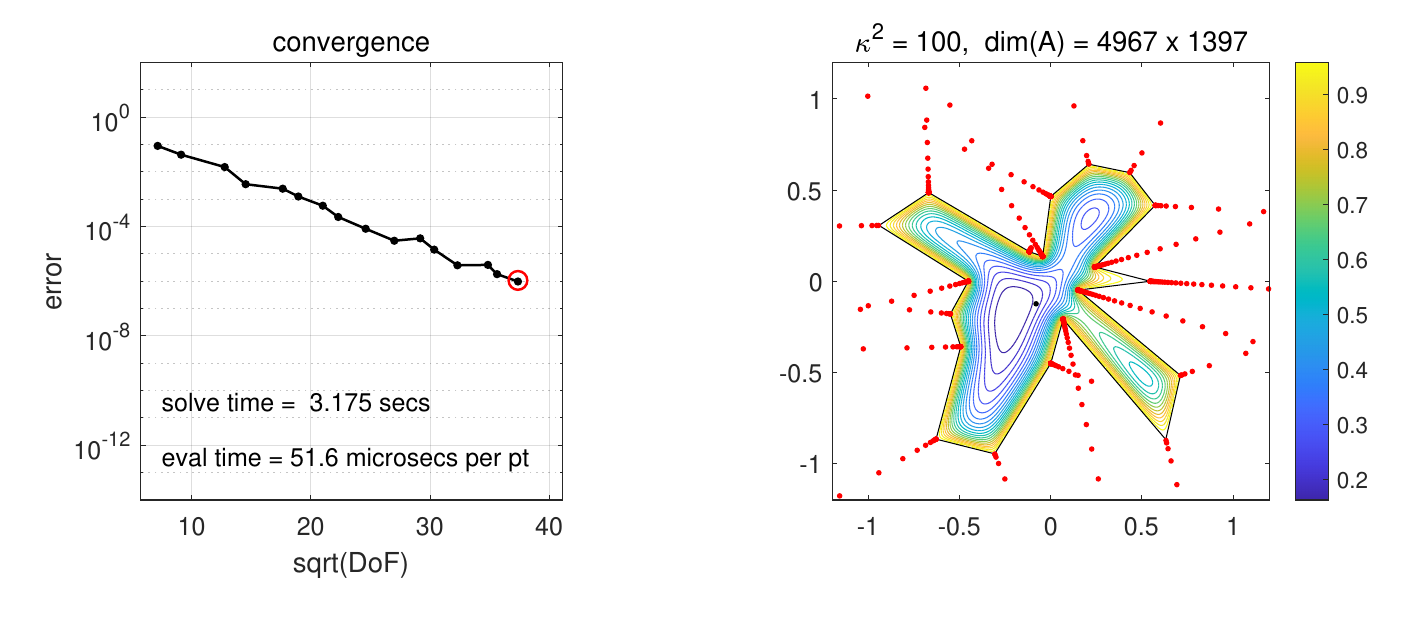}
    \vspace*{-10pt}
    \caption{Modified Helmholtz tests on a pentagon ($\kappa^2=1000$), a snowflake
    ($\kappa^2=10^4$), and a random 20-sided polygon ($\kappa^2=100$).  Straight curves
    in the semilog plots against $\sqrt N$ indicate root-exponential decay.  Red
    markers show independent boundary validation residuals rather than certified
    continuous-boundary errors.}
    \label{fig:poi-type-more}
\end{figure}

\cref{fig:poi-type-more} presents a pentagon with $\kappa^2=1000$, a snowflake with $\kappa^2=10^4$, and a random 20-sided polygon with $\kappa^2=100$.  
The solutions are nonnegative and exhibit increasingly thin boundary layers as $\kappa$ grows.  
The semilog plots against $\sqrt N$ again display approximately straight residual curves.
The pentagon and snowflake examples reach a validation residual of approximately $10^{-6}$ in less than one second, while the 20-sided polygon requires about eight seconds under the reported implementation.

\subsection{Asymptotic behaviors of coefficients as $\kappa \to 0$}\label{subsec:alpha-zero}

\begin{figure}[!h]
    \centering
    \subfigure[$\|a_n\|+\|b_n\|$ versus $\kappa$]{\includegraphics[width=0.45\linewidth]{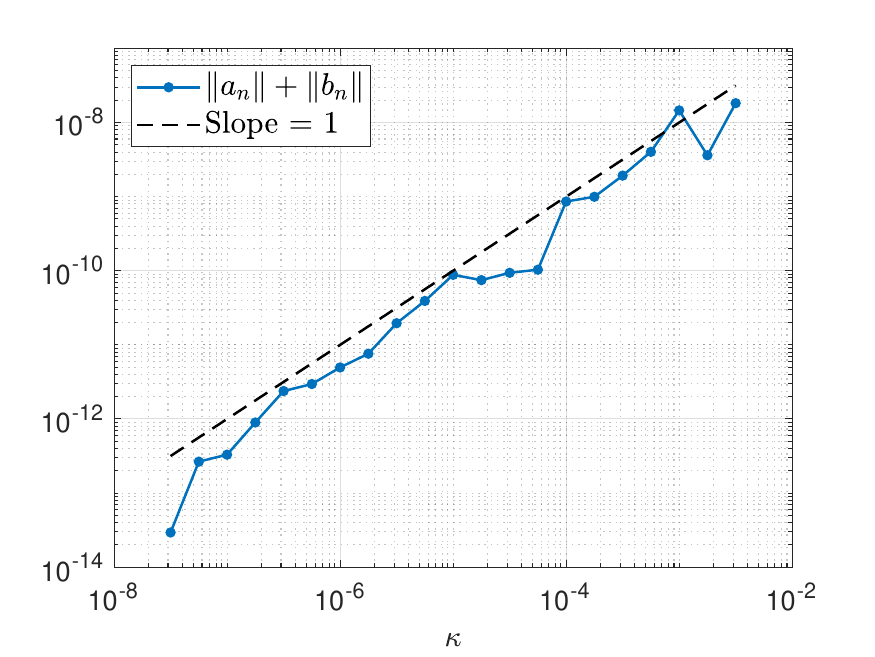}} 
    \subfigure[$c_m$ versus $m$]{\includegraphics[width=0.45\linewidth]{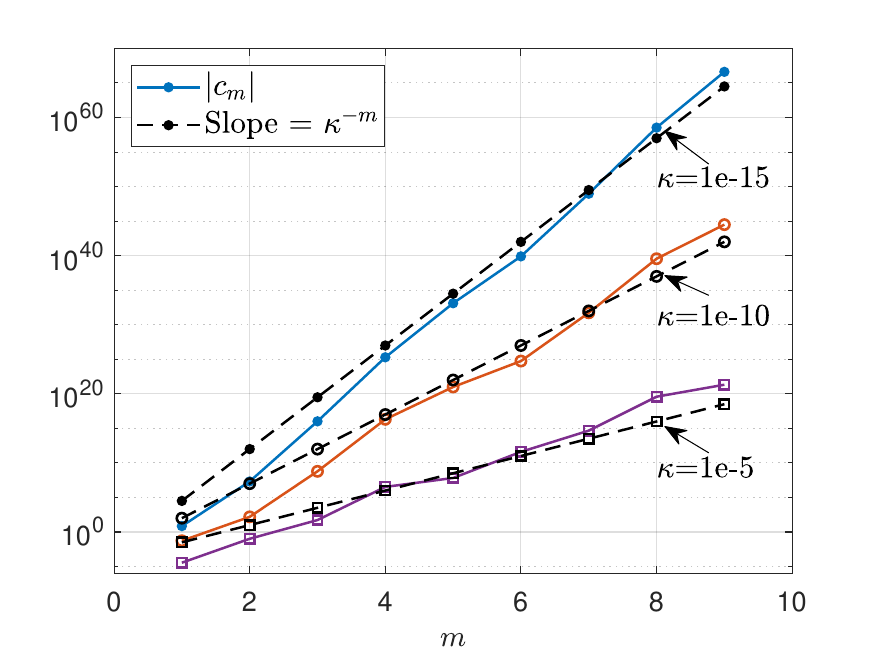}}
    \vspace*{-10pt}
    \caption{Coefficient behaviors in the limit $\kappa^2\to0$, compared with the
    scalings in~\eqref{eq:an-bn-alp0}.}
    \label{fig:poitype-alp0}
\end{figure}

Formally, the Laplace equation can be seen as the limit of the modified Helmholtz equation \eqref{eq:prob10} as $\kappa \to 0$. 
Therefore, it is expected that the approximation scheme \eqref{eq:poi-app} converges to \eqref{eq:rn-app} as $\kappa \to 0$. 
According to the asymptotic behaviors of the modified Bessel functions in \eqref{eq:asymp-IK}, we have
for an exterior pole $z_n$, 
\begin{equation*}%
  K_1(\kappa\, r_n)\re^{-\ri\theta_n} \sim \frac{1}{\kappa\,(z-z_n)}, \quad
  K_1(\kappa\, r_n)\re^{\ri\theta_n} \sim \frac{1}{\kappa\,(\overline z-\overline z_n)},
\end{equation*}
and for $m\geq1$, 
\begin{align*}%
 I_m(\kappa\, r_*)\re^{\ri m\theta_*} \sim \frac{\kappa^m}{2^m m!}(z-z_*)^m,\quad
 I_m(\kappa\, r_*)\re^{-\ri m\theta_*} \sim \frac{\kappa^m}{2^m m!}(\overline z-\overline z_*)^m,
 \nonumber
\end{align*}
while $I_0(\kappa r_*)\to1$.  
Hence the modified-Bessel ansatz tends, after coefficient
rescaling, to analytic and anti-analytic rational/polynomial harmonic fields.  The
corresponding coefficient scalings are
\begin{equation}\label{eq:an-bn-alp0}
  \|a_n\|+\|b_n\|=O(\kappa),\qquad
  |c_{m}|,\; |d_m|=O(\kappa^{-m}). %
\end{equation}
In \cref{fig:poitype-alp0}, we compute the coefficients $\|a_n\|+\|b_n\|$ and $|c_m|$ in \eqref{eq:poi-app} versus $\kappa$ when simulated the modified Helmholtz equation defined on a square domain with a side length of $1$.
It is obvious that the asymptotic behaviors of the coefficients coincide with the estimates in \eqref{eq:an-bn-alp0}.

\section{Lightning algorithm for inhomogeneous modified Helmholtz equations}\label{sec:inhom-poi}
In this section, we turn our attention to inhomogeneous modified Helmholtz equations.
The central task is the efficient computation of a particular solution, for which we introduce two classes of embedded domain methods.
With these solvers in place, we then investigate their applications to more general PDEs in the next section.

\begin{figure}[!t]
    \centering
    \subfigure[Circular embedding]{\includegraphics[width=0.45\linewidth]{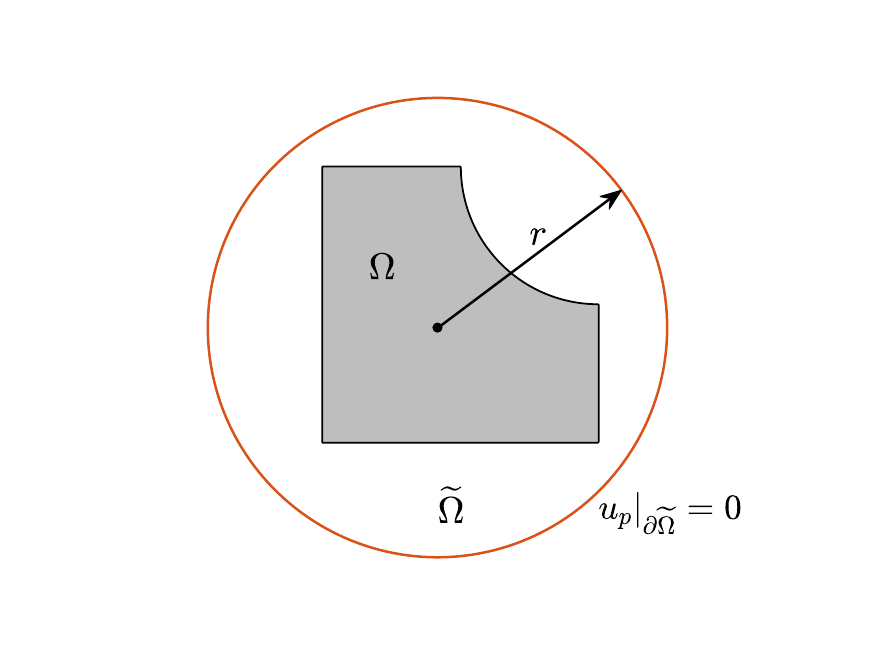}}\quad
    \subfigure[Fourier continuation]{\includegraphics[width=0.45\linewidth]{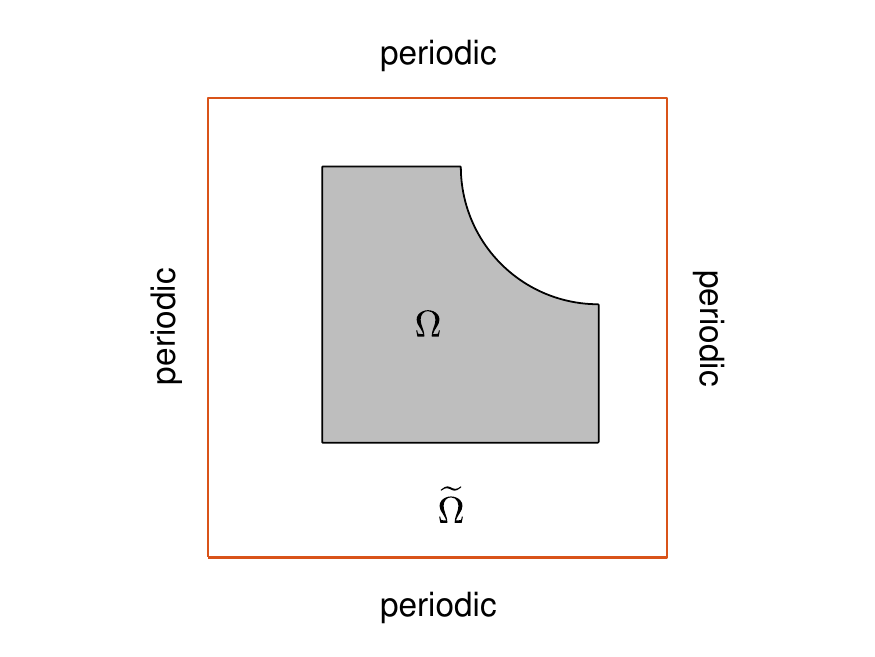}}
    \vspace*{-10pt}
    \caption{Illustration of the embedded domain method.}
    \label{fig:embed-FC}
\end{figure}

With an efficient algorithm for solving the homogeneous problem \eqref{eq:prob10}, one can compute the solution to the modified Helmholtz equation \eqref{eq:prob1} in $\Omega$ once a particular solution of \eqref{eq:sol-part} is efficiently computed. 
Actually, considerable specialist literature exists on methods for solving particular solutions.
Here, we introduce two standard techniques in the embedded domain methods for computing the particular solution \eqref{eq:sol-part}: (i) {\it embedding method}, and (ii) {\it Fourier continuation method}.
The former aims to embed the irregular domain $\Omega$ into a regular one, such as a disk or rectangle, where various efficient numerical methods can be applied, provided that appropriate boundary conditions are imposed.
While the Fourier continuation method allows for the use of the fast Fourier transform (FFT) by extending the target function to a rectangular domain periodically. 
The primary challenge in both methods is smoothly extending a function (typically given at a set of distinct points) to a regular domain.
Fortunately, significant progress has been made in addressing this challenge at a reasonable computational cost. Researchers such as Bruno, Lyon, and Huybrechs \cite{bruno2026two,lyon2011fast, bruno2022two, matthysen2018function} have successfully tackled this issue.
We will explore these techniques, in combination with our proposed algorithm, in future works. 
In the meantime, we assume that a smooth extension exists with suitable smoothness for the particular solution.
In \cref{fig:embed-FC}, we illustrate the main ideas behind the two techniques, showcasing how each method handles the extension of functions to regular domains.

\subsection{Embedding method}\label{subsec:embed}
As illustrated in \cref{fig:embed-FC} (a), the embedding approach for the modified Helmholtz equation in \eqref{eq:sol-part} seeks a solution $\tilde u_p$ that satisfies
\begin{equation}\label{eq:mod-embed}
    -\Delta \tilde u_p + \kappa^2 \tilde u_p = \tilde f,\quad \text{in\;\;} \widetilde \Omega,
    \quad \tilde u_p|_{\partial \widetilde \Omega} = 0.
\end{equation}
Here, $\tilde f$ denotes the extension of the function $f$ onto the disk domain $\widetilde \Omega$, with the property $\tilde f|_{\Omega} = f$, and similarly, $\tilde u_p$ is the extension of $u_p$ onto $\widetilde \Omega$.
It is well-known that the problem in \eqref{eq:mod-embed} is well-posed, and that $\tilde u_p \in H^{m+2}(\widetilde \Omega)$ provided $f \in H^m(\widetilde \Omega)$ for some integer $m \ge 0$ (see, e.g., \cite{evans2022partial}).

To obtain an accurate approximation to the solution of \eqref{eq:mod-embed}, we employ the spectral Galerkin method as presented in \cite{Shen1997SISC}. First, we rewrite the equation \eqref{eq:mod-embed} in polar coordinates $(r,\theta)$,
\begin{equation}\label{eq:mod-embed2}
\begin{split}
    &-\frac{1}{r} \partial_r(r \partial_r \tilde u_p) - \frac{1}{r^2} \partial_\theta^2 \tilde u_p + \kappa^2 \tilde u_p = \tilde f,\quad (r,\theta) \in (0,1) \times [0,2\pi), \\
    &\tilde u_p(1, \theta) = 0, \quad \theta \in [0, 2\pi), \quad \tilde u_p \text{ periodic in $\theta$.}
\end{split}
\end{equation}
Here, we assume that $\Omega$ belongs to the unit disk $\widetilde \Omega$ for simplicity.
The Fourier-spectral approximation to the solution $\tilde u_p$ leads to the following system of equations
\begin{equation}\label{eq:mod-embed-rt}
\begin{split}
    &-\partial_r^2 \tilde u_p^{im} - \frac1r \partial_r \tilde u_p^{im} + \Big( \frac{m^2}{r^2} + \kappa^2 \Big) \tilde u_p^{im} = \tilde f^{im}(r),\quad 0 \le r \le 1,\;\; i=1,2, \\
    &\tilde u_{p}^{1m}(0) = \tilde u_{p}^{2m}(0) = 0 \text{ if } m \neq 0,\;\;
    \tilde u_{p}^{1m}(1) = \tilde u_{p}^{2m}(1) = 0.
\end{split}
\end{equation}
Here, $\{\tilde u_{p}^{1m}, \tilde u_{p}^{2m}\}$ represent the Fourier coefficients of $\tilde u_p$, i.e.,
\begin{equation*}
    \tilde u_p(r, \theta) \approx \sum_{m=0}^{N_\theta} \big( \tilde u_{p}^{1m}(r) \cos(m\theta) + \tilde u_{p}^{2m}(r) \sin(m\theta) \big)
\end{equation*}
where $N_\theta$ is a positive integer, and likewise for $\{\tilde f^{1m}, \tilde f^{2m}\}.$
The conditions $\tilde u_{p}^{1m}(0) = \tilde u_{p}^{2m}(0) = 0$ are the essential pole condition for ensuring the well-posedness of the problem. For more  details, we refer the reader to \cite{Shen1997SISC} .

The corresponding variational formulation of the problem \eqref{eq:mod-embed-rt} is to find $\tilde u_p^{im} \in X(m)$ such that
\begin{equation*}
    \int_0^1 \partial_r \tilde u_p^{im} \partial_r v r \rd r + \int_0^1 \Big( \frac{m^2}{r^2} + \kappa^2 \Big) \tilde u_p^{im} v r \rd r
    = \int_0^1 \tilde f^{im} v r \rd r, \quad \forall v \in X(m),\;\; i=1,2,
\end{equation*}
where $X(m)=H_0^1(0,1)$ if $m \neq 0$ and $X(0)=\{v \in H^1(0,1):\; v(1) = 0\}.$
Based on this variational formulation, one can derive the spectral-Galerkin scheme using appropriate basis functions, such as the Legendre or Chebyshev polynomials \cite{Shen1997SISC}.

\subsection{Fourier continuation method}
The core idea behind the Fourier continuation method applied to the modified Helmholtz equation in \eqref{eq:sol-part} involves periodically extending the function to a domain, such as the box $\widetilde \Omega=[0,2\pi)^2$, as depicted in \cref{fig:embed-FC} (b). 
This extension allows for the application of the Fourier spectral method, which benefits from the efficiency of the FFT. We omit the detailed derivation of the Fourier spectral method here and refer the reader to \cite{Trefethen2000Book,Shen2011Book} for a comprehensive discussion of its implementation and theoretical analysis.

\subsection{Numerical experiments}
\begin{figure}[!t]
    \centering
    \includegraphics[width=0.8\linewidth]{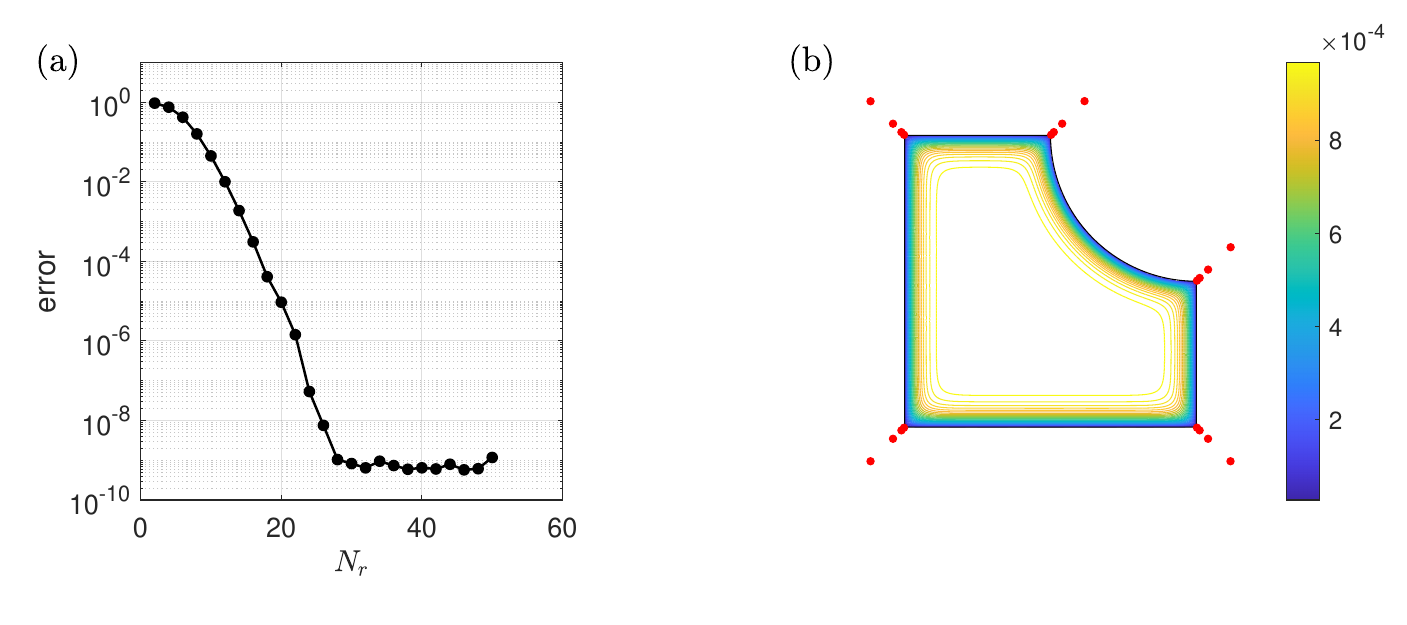}
    \includegraphics[width=0.8\linewidth]{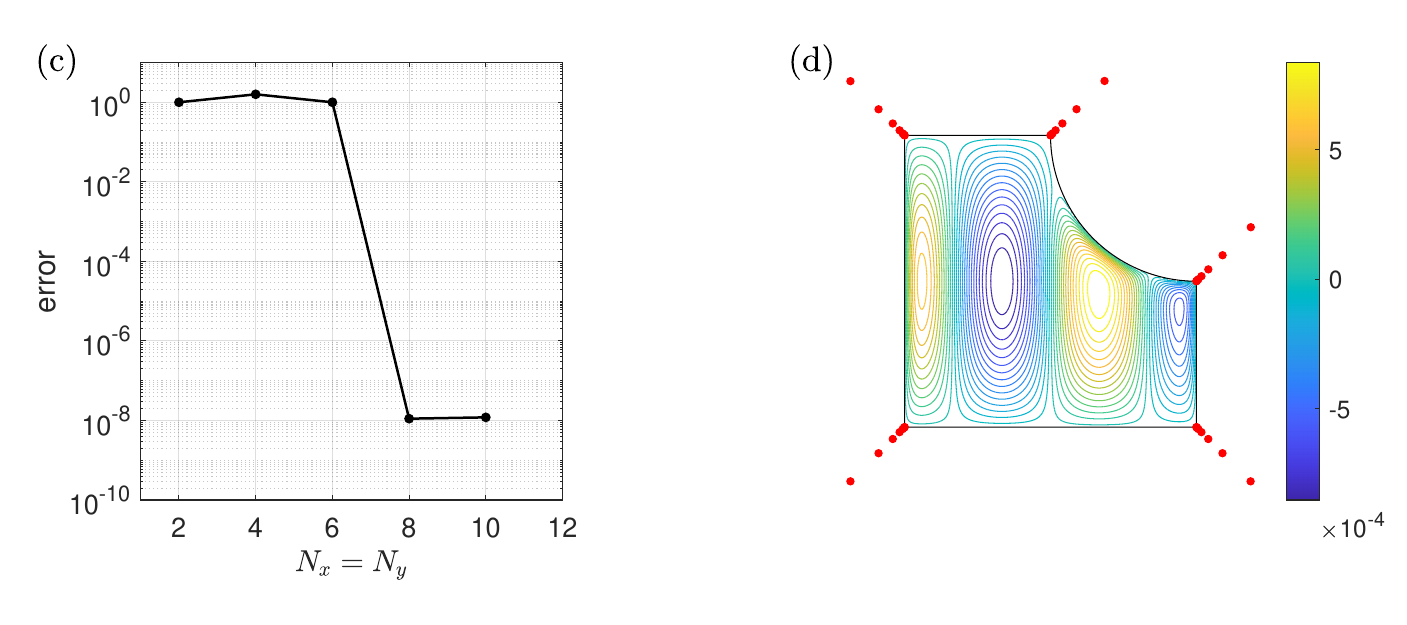}
    \vspace*{-10pt}
    \caption{Numerical results for the inhomogeneous modified Helmholtz equation 
    in \eqref{eq:sol-part}. (a)-(b): the embedding method; (c)-(d): the Fourier continuation method.}
    \label{fig:inhom-embedding}
\end{figure}

To assess the efficiency and accuracy of the methods, we set $\kappa^2=1000$ and consider an L-shaped domain with one curved edge, as illustrated in \cref{fig:embed-FC}.
In addition, the following settings are chosen for the embedding method and the Fourier continuation method, respectively:
\begin{itemize}
    \item [(i)] $f(x,y)=1$ and $h=0$ on $\partial \Omega$, $N_\theta=100$ in $\theta$ and Legendre spectral Galerkin method in $r$ with degree $N_r$;
    \item [(ii)] $f(x,y) = \sin(3\pi x) \cos(\pi y)$ and Fourier modes $N_x=N_y$ in $x$ and $y$.
\end{itemize}
Numerical results are presented in \cref{fig:inhom-embedding}. The first row (a)-(b) shows the maximum-norm errors as a function of $N_r$ and the contour plot of the numerical solution, respectively, for the embedding method. The second row (c)-(d) displays the corresponding results for the Fourier continuation method.
The red dots indicate the location of poles.
As evident, spectral accuracy is achieved even in the presence of singularities near corners.

\section{Applications}\label{sec:app}
We next use the inhomogeneous solver as the elliptic kernel in time-dependent computations.  For fixed time step, diffusion coefficient, reaction parameter, and geometry, the boundary least-squares matrix is unchanged from one time level to the next.  Its column scaling and rank-revealing factorization can therefore be computed once and reused for all subsequent right-hand sides.  This separation between geometry-dependent setup and per-step solution cost is an important feature of the method.

For large reaction parameters, boundary layers may arise in addition to corner singularities. Accordingly, the pole sets used below consist of the standard corner clusters together with additional exterior clusters placed near selected boundary-layer locations. In our simulations, these locations are treated as `pseudo-singular' points along the boundary, around which additional poles are clustered to better resolve the boundary layers.
Adaptive pole selection by AAA-type procedures is another possibility; see \cite{Nakatsukasa2018SISC,Costa2023aaa,xue2024stokes}.

\subsection{RBF frame for function expansion}\label{subsec-rbf}

\begin{figure}[!h]
    \centering
    \includegraphics[width=0.5\linewidth]{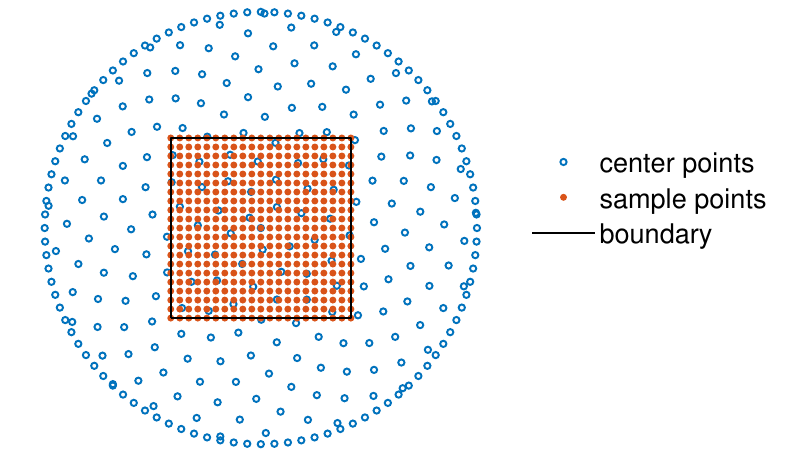}
    \caption{Illustration of the center points and sample points for the
    RBF frame.}
    \label{fig:rbf-fram}
\end{figure}

In the simulations of this section, we employ the RBF frame as the function expansion basis \cite{adcock2022stable}. The physical domain is fixed as $\bs x=(x,y)\in\Omega:=[-\frac12,\frac12]^2$.
Given a set of sample data $\{\bs x_j,f_j\}_{j=1}^{N_s}$ and a collection of center points $\{\bs y_k\}_{k=1}^{N_c}$ ($N_c \gg N_s$), we approximate the function by the RBF expansion
\begin{align*}
    \tilde f(\bs x) = \sum_{k=1}^{N_c} s_k \phi(\varepsilon,\bs x-\bs y_k),
\end{align*}
where the unknown coefficient vector $\bs s=(s_1,\ldots,s_{N_c})^\top\in\mathbb R^{N_c}$ is determined by the least-squares problem
\begin{equation*}
    \bs s = \operatorname*{arg\, min}_{\bs c\in\mathbb R^{N_c}} \|\bs\Phi\bs c-\bs f\|_2,
\end{equation*}
with the interpolation matrix $\bs\Phi=(\Phi_{jk}) = \bigl( \phi(\varepsilon,\bs x_j-\bs y_k) \bigr) \in\mathbb R^{N_s\times N_c},$ and $\bs f=(f_1,\ldots,f_{N_s})^\top.$
The distributions of the center points and sample points are illustrated in Figure~\ref{fig:rbf-fram}.

In the computation, we use Gaussian RBFs as basis functions, $\phi_k(\varepsilon,\bs x)=\exp\bigl(-\varepsilon^2 \|\bs x-\bs y_k\|_2^2\bigr),$ and set $\varepsilon=0.1\sqrt{2N_c}.$
We select $3,000$ Fibonacci spiral points uniformly distributed in the unit disk as the center points and $6,400$ uniformly distributed points in $\Omega$ as the sample points.

\subsection{Applications to PDEs}
For general parabolic equations, we consider the following model problem
\begin{equation*}
    \partial_t u = \nu \Delta u + \mathcal N(u,t), \quad \text{in }\;\; \Omega, \quad t \in (0, T],
\end{equation*}
with suitable initial and boundary conditions, we uniformly partition the time interval $[0,T]$ into $N_t$ subintervals with time step $\tau=T/N_t$, and compute the first-step approximation $u^1$ using the implicit-explicit Euler scheme
\begin{equation}\label{eq:app-u1}
    - \Delta u^1 + \frac{1}{\tau \nu} u^1  = \frac{1}{\tau \nu } \big( u^0 + \tau \mathcal N(u^0, t_0)\big)=: g^1,
\end{equation}
and compute the approximations $u^{n+1}$ for $n\geq1$ subsequently using the implicit-explicit BDF2 scheme
\begin{equation}\label{eq:app-un}
    - \Delta u^{n+1} + \frac{3}{2\tau \nu} u^{n+1}  = \frac{1}{2\tau \nu} \big( 4u^n-u^{n-1} + 2\tau \mathcal N(u^n, t_n)\big)=:g^{n+1}, 
\end{equation}
for $1 \le n \le N_t-1.$
For inhomogeneous Helmholtz equations \eqref{eq:app-u1}-\eqref{eq:app-un},  we can compute the known $u^n$ sequentially by the method proposed in previous sections.
We summarize our method in \cref{alg:app} and consider its applications, such as the diffusion equation and the viscous Burgers' equation, at the end of this section.

\begin{algorithm}[!t]
\caption{Lightning algorithm for PDEs on complex domains}\label{alg:app}
\begin{algorithmic}[1]
\Require Domain $\Omega$, parameters $\nu$, final time $T$, positive integer $N_t$, nonlinear term $\mathcal N$, initial and boundary conditions.
\Ensure Numerical solution $u^{N_t}$.
\State Set $u^{n+1} = u_p + v$ with $u_p$ and $v$ satisfying the equation \eqref{eq:sol-part} and \eqref{eq:prob10}, respectively.
\State For \eqref{eq:app-u1}, extend $g^1$ by the RBF method proposed in  \cref{subsec-rbf}, compute $u_p$ by spectral Galerkin methods, and compute $v$ by the lightning algorithm.
\For{$n=1,\ldots,N_t-1$}
 \State Given $u^n$ and $u^{n-1}$, extend $g^{n+1}$  in \eqref{eq:app-un} by the RBF extension method, and compute $u_p$ and $v$, sequentially.
\EndFor
\end{algorithmic}
\end{algorithm}

\subsubsection{Diffusion equation}
In our first example, we focus on the diffusion equation:
\begin{equation}\label{eq:diff-model}
\begin{dcases}
  \partial_tu-\nu\Delta u=f,&\text{in }\Omega\times(0,T],\\
  u(\cdot,0)=u_0,&\text{in }\overline\Omega,\\
  u=h,&\text{on }\partial\Omega\times[0,T].
\end{dcases}
\end{equation}
Let $t_n=n\tau$, $N_t\tau=T$, and let $u^n$ approximate $u(t_n)$. Then applying the implicit-explicit BDF2 scheme gives
\begin{equation}\label{eq:diff-bdf2}
\begin{dcases}
  -\Delta u^{n+1}+\dfrac{3}{2\tau\nu}u^{n+1}=g^{n+1},&\text{in }\Omega,\\
  u^{n+1}=h^{n+1},&\text{on }\partial\Omega,
\end{dcases}
\end{equation}
where
\begin{equation*}%
  g^{n+1}=\frac1\nu\Big(f^{n+1}+\frac{4u^n-u^{n-1}}{2\tau}\Big), \quad n=1,\ldots,N_t-1.
\end{equation*}
The first value $u^1$ is computed by the backward Euler scheme.

\begin{figure}[H]
    \centering
    \includegraphics[width=0.33\linewidth]{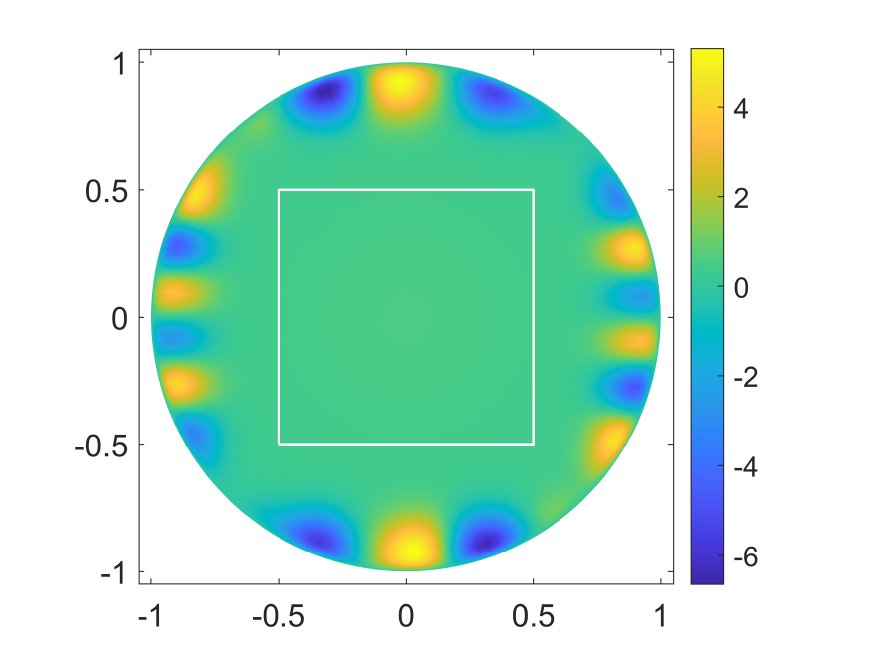} \hspace*{-5pt}
    \includegraphics[width=0.33\linewidth]{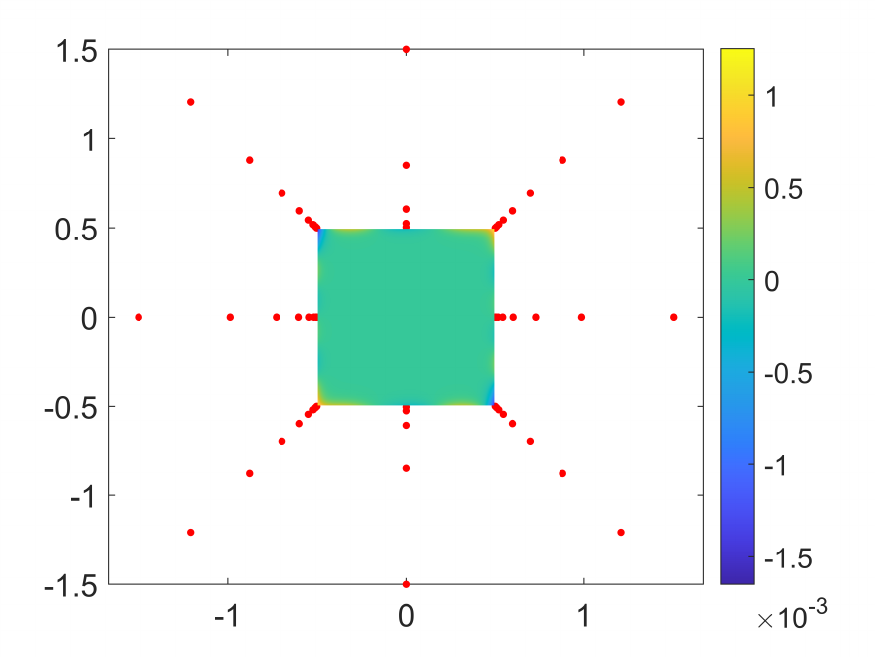} \hspace*{-5pt}
    \includegraphics[width=0.33\linewidth]{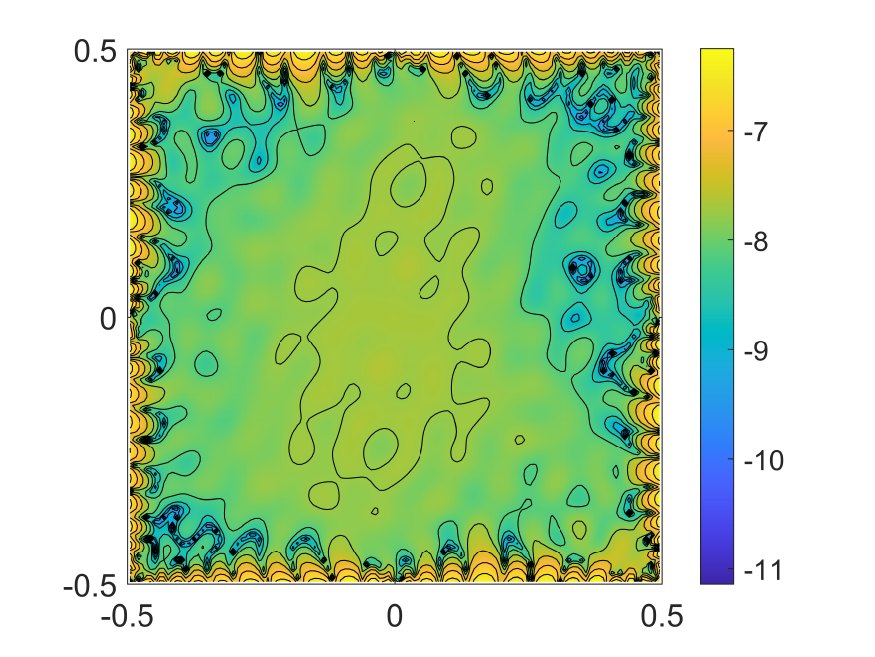} \\
    \subfigure[$u_p$]{\includegraphics[width=0.33\linewidth]{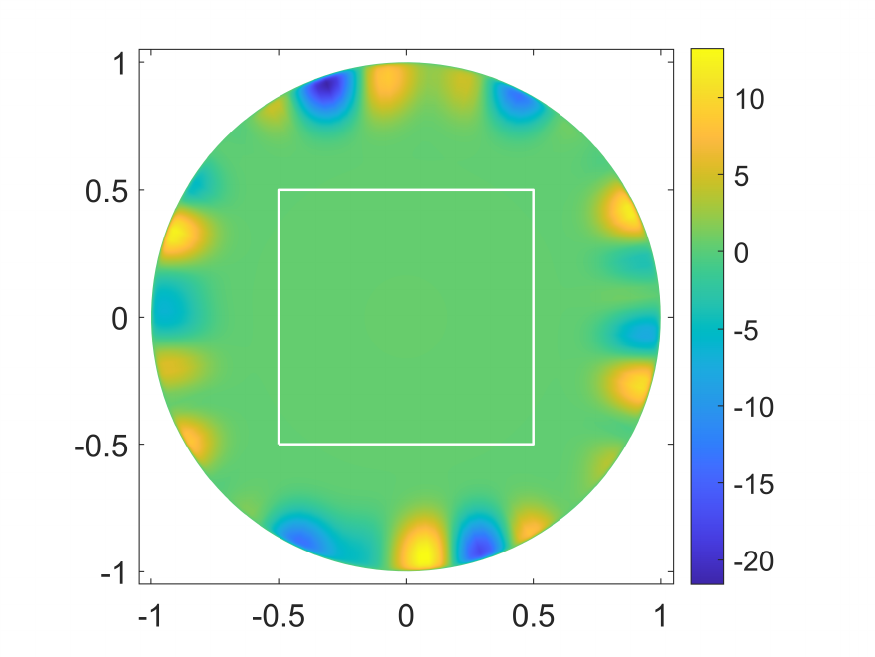}} \hspace*{-5pt}
    \subfigure[$v$]{\includegraphics[width=0.33\linewidth]{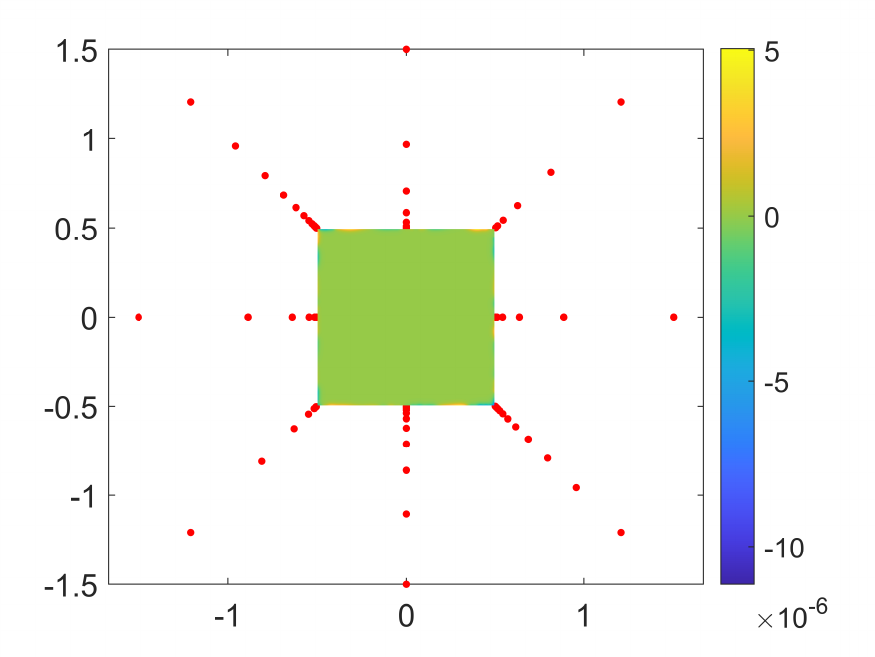}} \hspace*{-5pt}
    \subfigure[$\log_{10}$(Error)]{\includegraphics[width=0.33\linewidth]{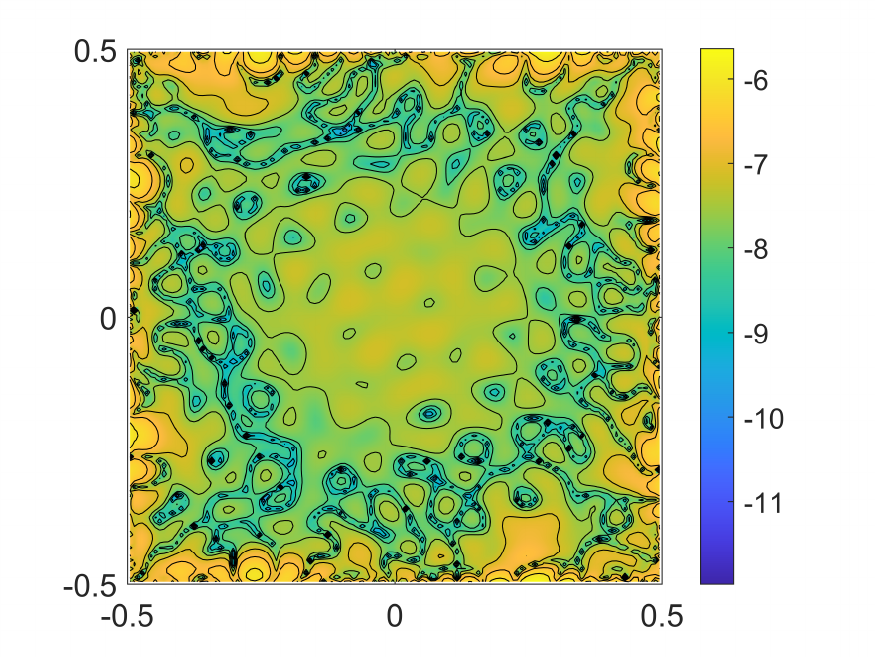}}
    \vspace*{-5pt}
    \caption{Diffusion problem~\eqref{eq:diff-model} with $\nu=1,$ $f=1$, $h=0$, and $u_0=0$ on a square domain. (a) Particular solution $u_p$ by RBF approximation; (b) Lightning algorithm for $v$ with tolerance $10^{-6}$; (c) Pointwise error in logarithmic scale.
    First row: $\nu=1$; Last row: $\nu=0.1.$}
    \label{fig:diff-err}
\end{figure}

To validate the accuracy of the proposed algorithm, we take $T=1$, $\tau=10^{-3}$, and choose $f$, $h$, and $u_0$ so that the exact solution is given by
\begin{equation}\label{eq:diff-init}
  u(x,y,t)=\frac{\cos t}{1+x^2+y^2}, \qquad (x,y)\in \Omega.
\end{equation}
The particular-solution right-hand side $g^{n+1}$ is extended to the unit disk by an RBF construction and solved by the polar spectral method.  
In \cref{fig:diff-err}, we show the profiles of the particular solution $u_p$ obtained by the RBF frame approximation, the correction $v$ computed by the lightning algorithm, and the pointwise errors on a logarithmic scale for $\nu=1$ (top row) and $\nu=0.1$ (bottom row), respectively.
In these computations, one pseudo-pole is added along each side of the square to improve the accuracy of the lightning approximation.
We observe that the largest pointwise errors are concentrated near the boundary. Moreover, as $\nu$ decreases, the large-amplitude features of the particular solution become increasingly localized near the boundary, and the lightning approximation requires more poles to resolve the resulting boundary behavior accurately.

\begin{figure}[!t]
    \centering
    \includegraphics[width=0.33\linewidth]{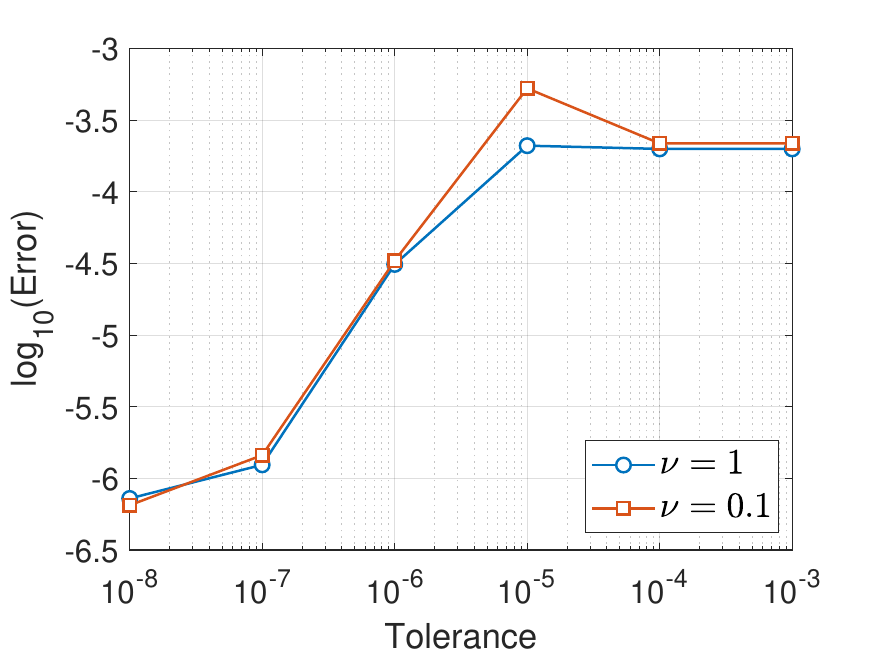} \hspace*{-5pt}
    \includegraphics[width=0.33\linewidth]{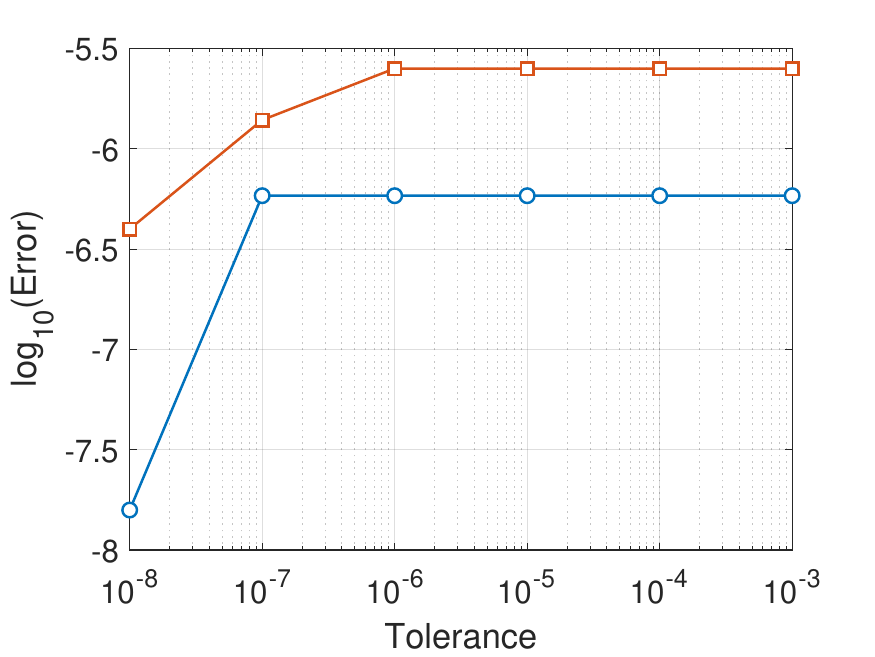} \hspace*{-5pt}
    \includegraphics[width=0.33\linewidth]{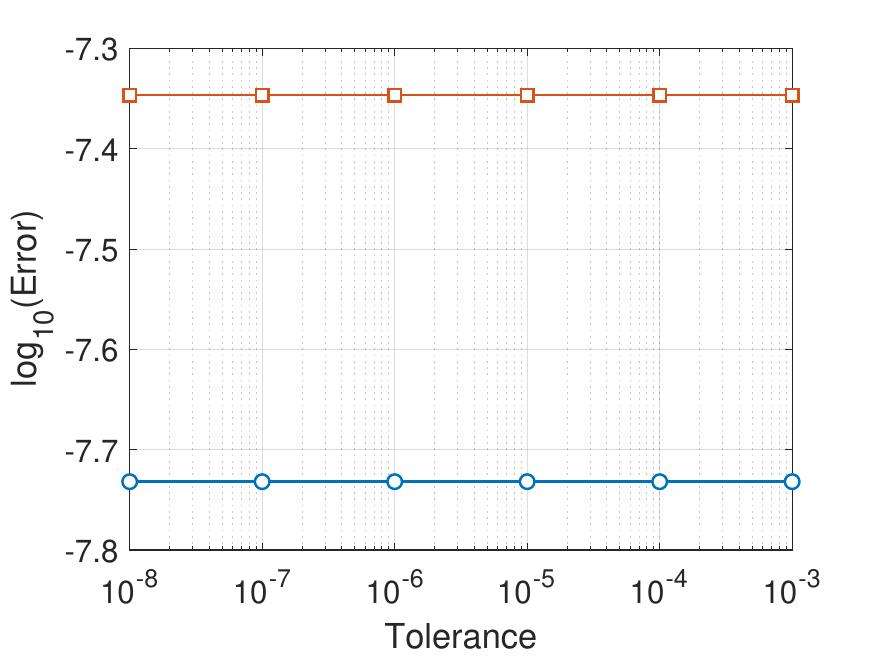}
\begin{tikzpicture}[scale=1]
\draw[semithick] (-5,0) node[left=0.1cm, font=\footnotesize] {(a)} -- (-3,0);
\filldraw[black] (-4,0) circle (1.5pt);

\draw[semithick] (-1,0) node[left=0.1cm, font=\footnotesize] {(b)} -- (1,0);
\filldraw[black] (-1/2,0) circle (1.5pt);
\filldraw[black] (0,0) circle (1.5pt);
\filldraw[black] (1/2,0) circle (1.5pt);

\draw[semithick] (3,0) node[left=0.1cm, font=\footnotesize] {(c)} -- (5,0);
\filldraw[black] (3+1/3,0) circle (1.5pt);
\filldraw[black] (3+2/3,0) circle (1.5pt);
\filldraw[black] (4,0) circle (1.5pt);
\filldraw[black] (4+1/3,0) circle (1.5pt);
\filldraw[black] (4+2/3,0) circle (1.5pt);
\end{tikzpicture}
\vspace*{-5pt}
    \caption{Maximum errors versus tolerances in the lightning algorithm under different number of pseudo-poles along each side. From left to right, the numbers of pseudo-poles are 1,3,5 on each side, respectively.}
    \label{fig:ill-pole}
\end{figure}

To further investigate the influence of the pseudo-poles, in \cref{fig:ill-pole} we report the maximum error at $T=1$ for different numbers of pseudo-poles placed along each side, with the tolerance of the lightning algorithm varying from $10^{-8}$ to $10^{-3}$, for both $\nu=1$ and $\nu=0.1$. The results show that increasing the number of pseudo-poles significantly reduces the sensitivity of the numerical error to the prescribed tolerance. In particular, with sufficiently many pseudo-poles, the error remains nearly unchanged over the entire range of tolerances considered, indicating improved robustness of the lightning algorithm.

\begin{figure}[!t]
    \centering
    \includegraphics[width=0.33\linewidth]{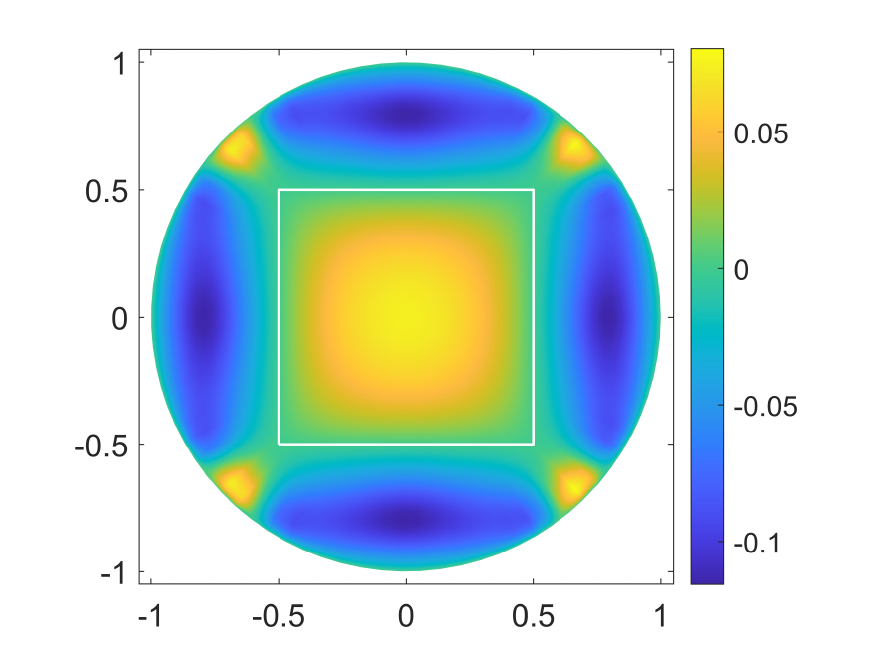} \hspace*{-5pt}
    \includegraphics[width=0.33\linewidth]{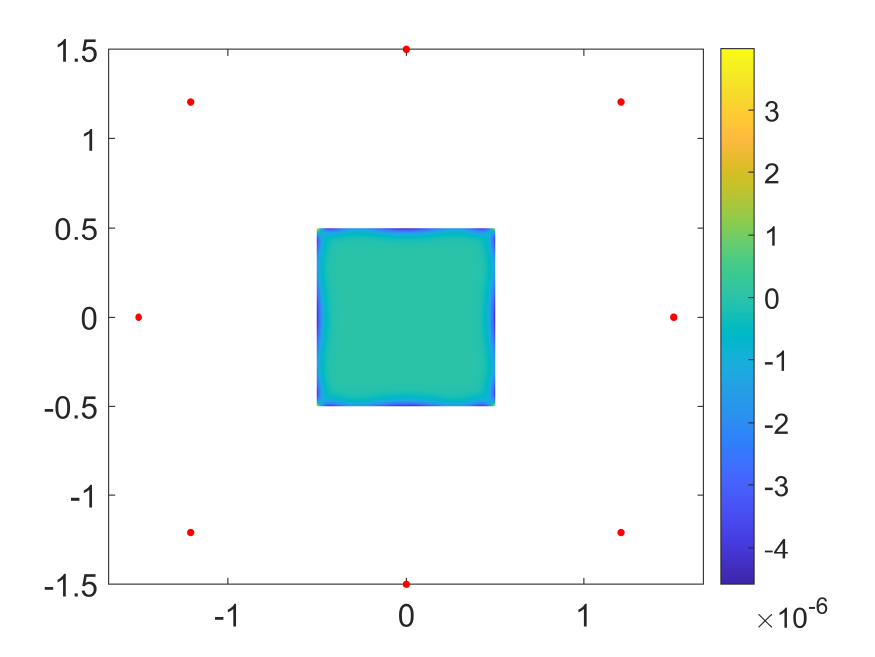} \hspace*{-5pt}
    \includegraphics[width=0.33\linewidth]{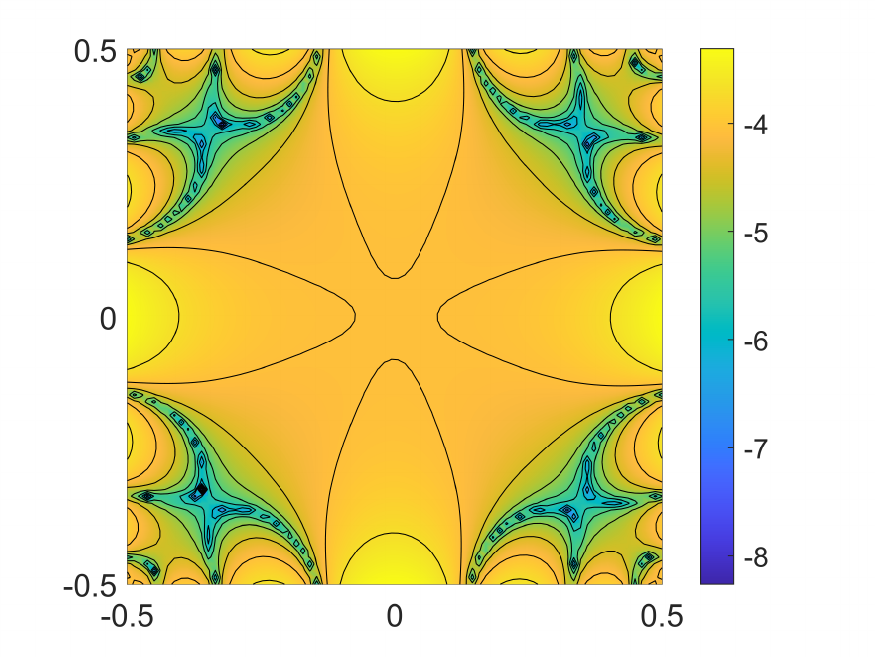} \\
    \subfigure[$u_p$]{\includegraphics[width=0.33\linewidth]{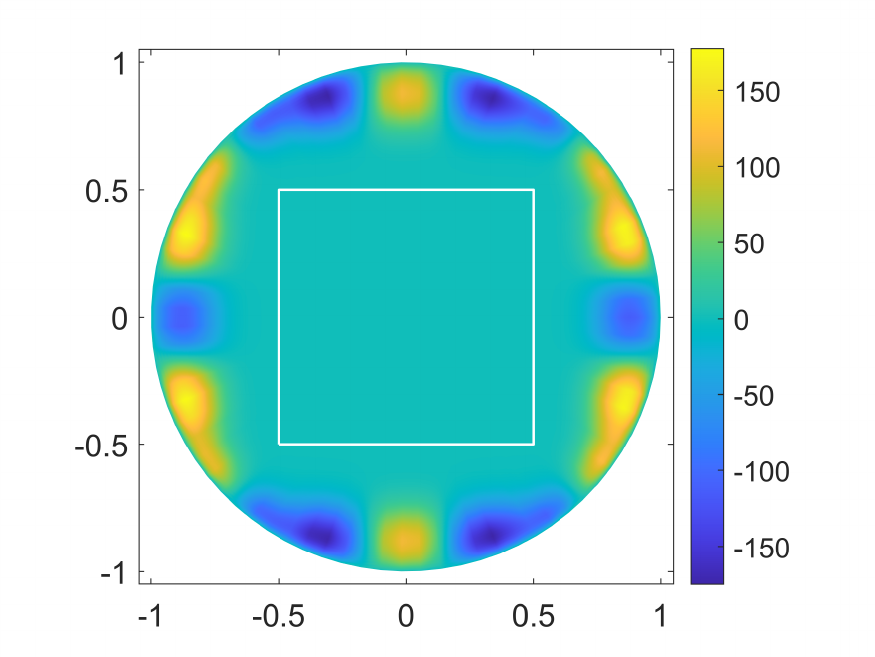}} \hspace*{-5pt}
    \subfigure[$v$]{\includegraphics[width=0.33\linewidth]{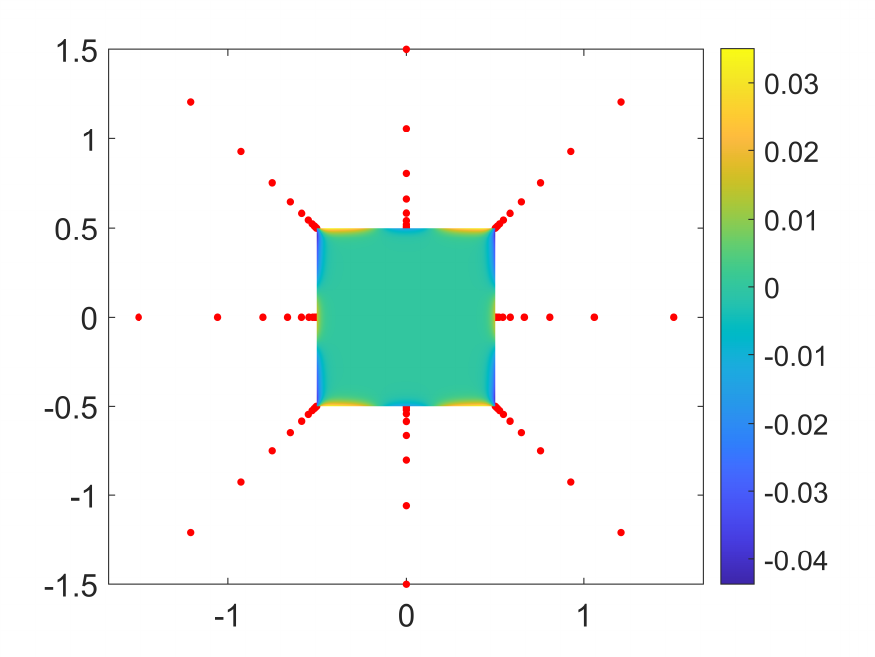}} \hspace*{-5pt}
    \subfigure[$\log_{10}$(Error)]{\includegraphics[width=0.33\linewidth]{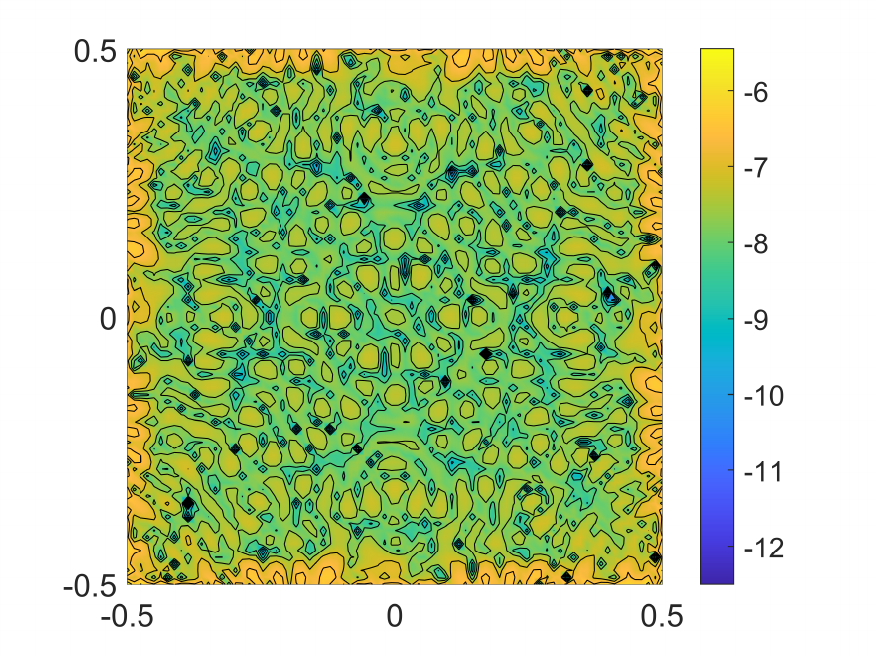}}
    \vspace*{-5pt}
    \caption{Similar to \cref{fig:diff-err} but with $\nu=1,$ $f=1$, $h=0$, and $u_0=0$ on a square domain and tolerance $10^{-3}$ (first row) and $10^{-6}$ (second row). The white line in (a) denotes the boundary $\partial \Omega.$}
    \label{fig:f1-comp2}
\end{figure}

\begin{figure}[!t]
  \centering
  \subfigure[Triangular domain]{\includegraphics[width=0.33\linewidth]{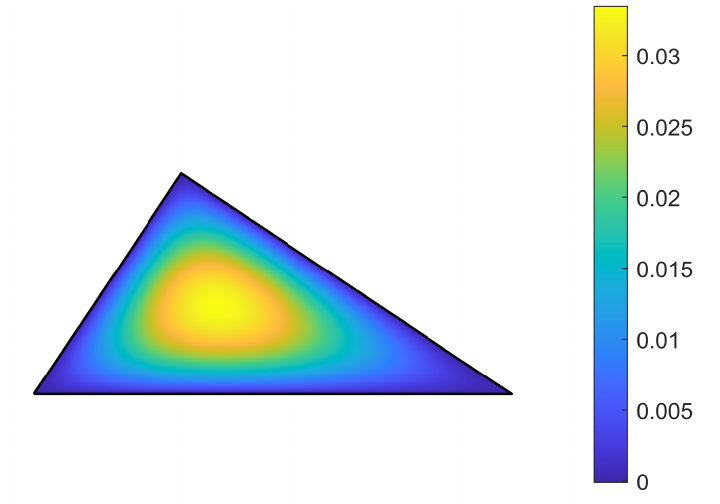}} \hspace*{-5pt}
  \subfigure[L-shaped domain]{\includegraphics[width=0.33\linewidth]{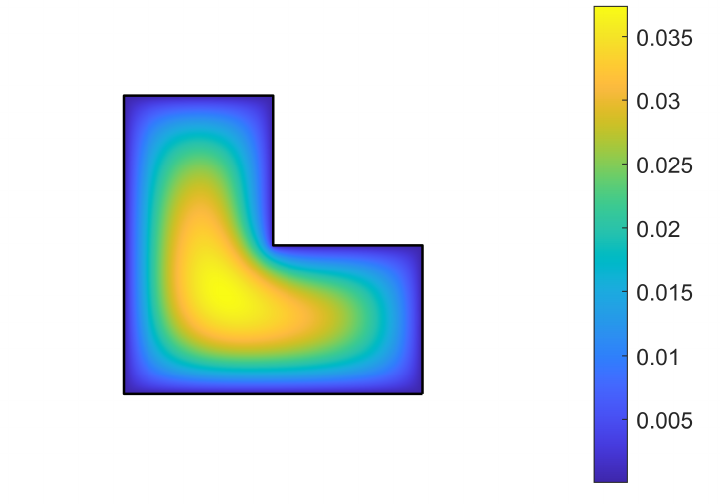}} \hspace*{-5pt}
  \subfigure[Snow-shaped domain]{\includegraphics[width=0.33\linewidth]{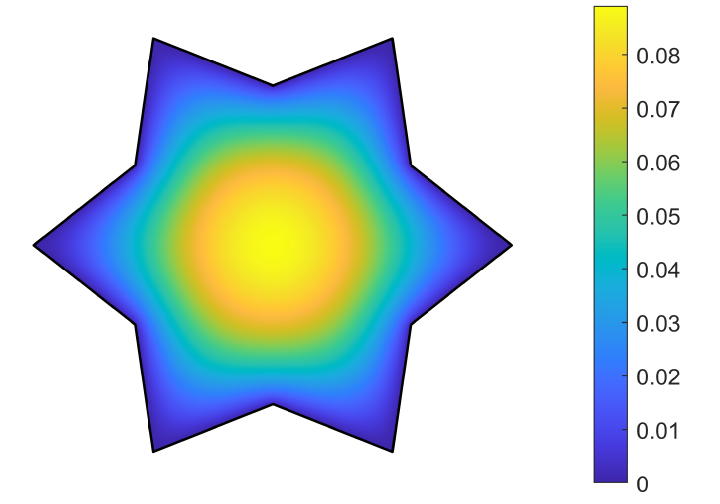}}
  \vspace*{-5pt}
  \caption{Similar to \cref{fig:f1-comp2} on three additional domains.}
  \label{fig:f1-more}
\end{figure}

We also consider the case $\nu=1$, $f=1$, $h=0$, and $u_0=0$, and use one pseudo-pole on each side (see \cref{fig:ill-pole} (a)). A high-resolution spectral--Galerkin solution is taken as the reference. In \cref{fig:f1-comp2}, we present the particular solution $u_p$ obtained by the RBF frame approximation, the correction $v$ computed by the lightning algorithm, and the corresponding pointwise errors for two different lightning tolerances, $10^{-3}$ and $10^{-6}$. The comparison shows that decreasing the tolerance improves the accuracy of the lightning approximation. 
In particular, the pointwise error is substantially reduced when the tolerance is decreased from $10^{-3}$ to $10^{-6}$.
In addition, we consider several further examples on different computational domains, as shown in \cref{fig:f1-more}.

\subsubsection{Viscous Burgers equation}
We next consider the viscous Burgers' equation
\begin{equation}\label{eq:burg}
\begin{dcases}
  \partial_tu-\nu\Delta u+\mu u(\partial_xu+\partial_yu)=0,
  &\text{in }\Omega\times(0,T],\\
  u(\cdot,0)=u_0,&\text{in }\overline\Omega,\\
  u=0,&\text{on }\partial\Omega\times[0,T].
\end{dcases}
\end{equation}
With $\widetilde u^n=2u^n-u^{n-1}$, the implicit-explicit BDF2 method is to solve
\begin{equation*}%
  -\Delta u^{n+1}+\dfrac{3}{2\tau\nu}u^{n+1} = g^{n+1} \quad \text{in }\Omega; \qquad
  u^{n+1}=0\quad \text{on }\partial\Omega,
\end{equation*}
with $g^{n+1}:= \tfrac1\nu \bigl(\frac{4u^n-u^{n-1}}{2\tau} - \mu\widetilde u^n (\partial_x \widetilde u^n + \partial_y \widetilde u^n)\bigr).$
Typically, the values of $u^1$ can be computed by the first-order implicit-explicit Euler method.

\begin{figure}[!t]
  \centering
  \subfigure[Initial function]{\includegraphics[width=0.33\linewidth]{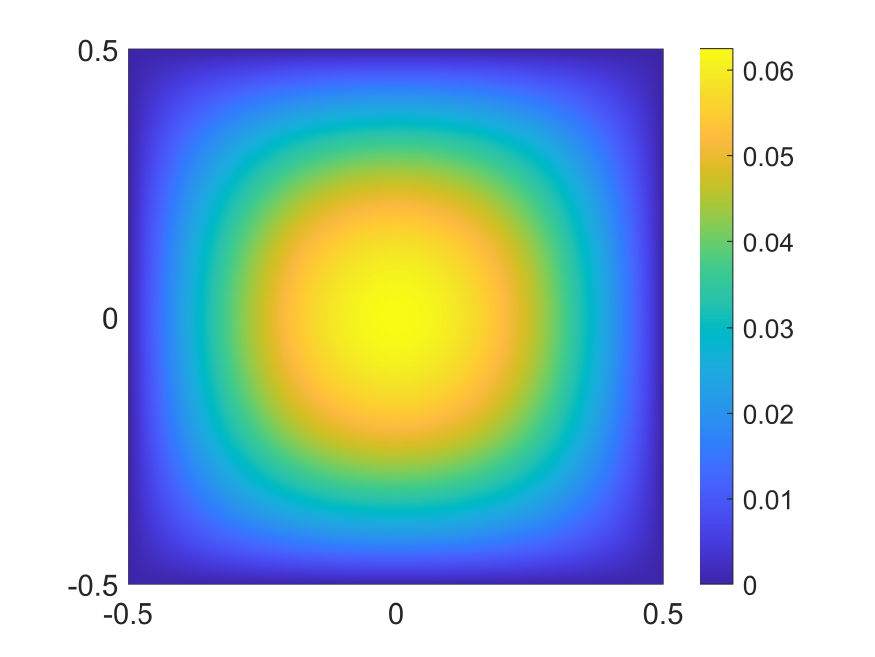}} \hspace*{-5pt}
  \subfigure[Numerical solution]{\includegraphics[width=0.33\linewidth]{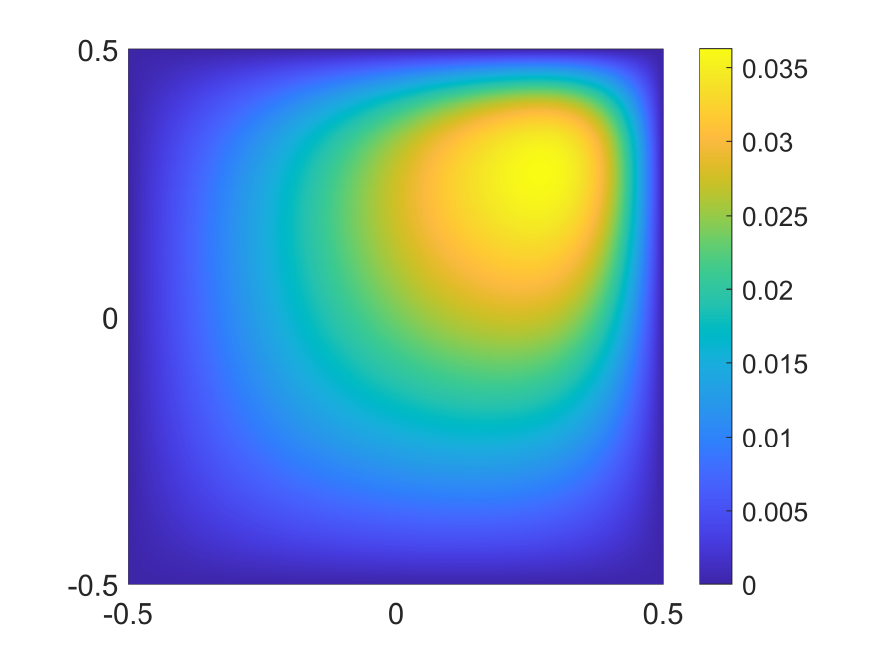}} \hspace*{-5pt}
  \subfigure[$\log_{10}$(Error)]{\includegraphics[width=0.33\linewidth]{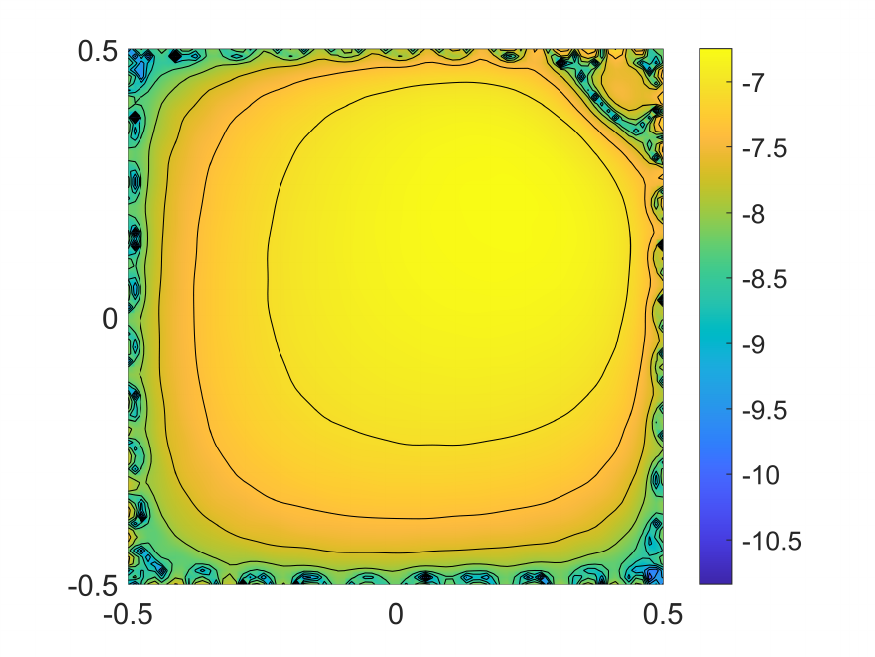}}
  \vspace*{-5pt}
  \caption{Square-domain viscous Burgers computation with $\nu=\tfrac{1}{200}$,
  $\mu=2$, $\tau=\tfrac{1}{400}$, $T=5$, and initial data~\eqref{eq:burg-init0}.}
  \label{fig:burg-err}
\end{figure}

For the square-domain test, we choose the initial profile as
\begin{equation}\label{eq:burg-init0}
  u_0(x,y)=(\tfrac12-x)^2(\tfrac12-y)^2, \quad (x,y)\in \Omega.
\end{equation}
A reference solution can be obtained by the spectral Galerkin method in space and a refined step-size in time.
In our simulation, we take $\nu=\tfrac{1}{200}, \mu=2$ and $\tau=\tfrac{1}{400}.$ Notably, in our lightning algorithm, we set the tolerance to be $10^{-7}$ and add nine pseudo-poles along each side.
In \cref{fig:burg-err}, we depict the initial profile given by \eqref{fig:burg-err}, the numerical solution at $T=5$ and the corresponding pointwise error.
It can be observed that the peak of the solution moves toward the boundaries resulting in a slightly sharper boundary layer, meanwhile our proposed algorithm gives an accurate approximation with the maximum pointwise error approximately $10^{-7}.$

Next, let us consider the approximation to \eqref{eq:burg} under different domains.
Particularly, we choose $\Omega$ to be a triangular domain and take the initial function as
\begin{equation}\label{eq:burg-init}
    u_0(x,y)=\sin^2(4\pi x) \sin^2(4\pi y) \chi(x,y),
\end{equation}
where $\chi(x,y)$ is a linear function with respect to $x$ and $y$ satisfying $\chi(x,y)=0$ on $\partial \Omega.$ To be more precise, assume the coordinates of a triangle is given by $(x_1, y_1), (x_2,y_2)$ and $(x_3,y_3),$ then we can formulate $\chi(x,y)$ as
\begin{align*}
    \chi(x,y) &= \big( (y_2-y_1)x - (x_2-x_1)y + x_2y_1 - x_1y_2 \big) \\
    &\qquad \times \big( (y_3-y_1)x - (x_3-x_1)y + x_3y_1 - x_1y_3 \big) \\
    &\qquad \times \big( (y_3-y_2)x - (x_3-x_2)y + x_3y_2 - x_2y_3 \big).
\end{align*}

\begin{figure}[!t]
  \centering
  \subfigure[$t=0$]{\includegraphics[width=0.33\linewidth]{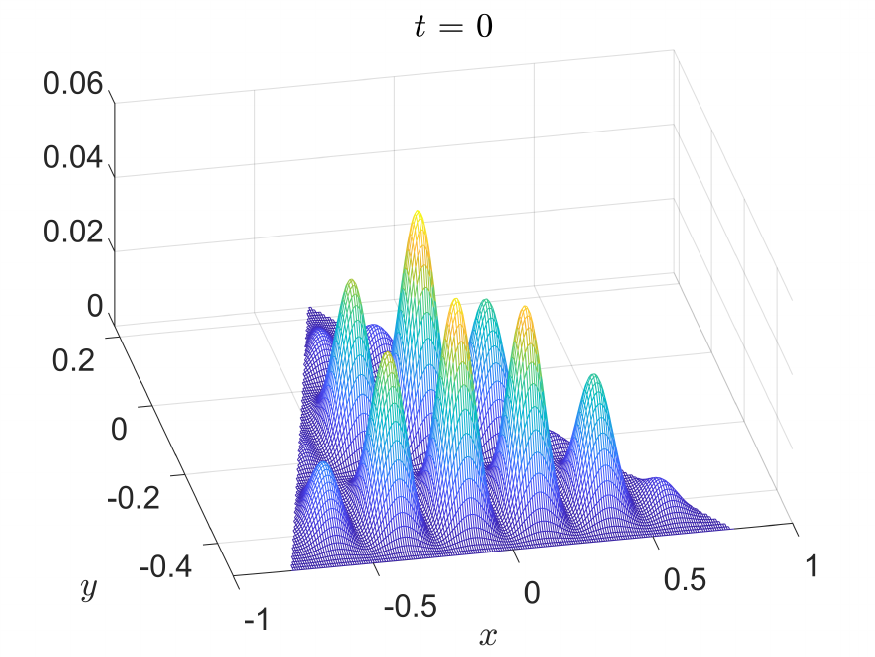}}\hspace*{-5pt}
  \subfigure[$t=3$]{\includegraphics[width=0.33\linewidth]{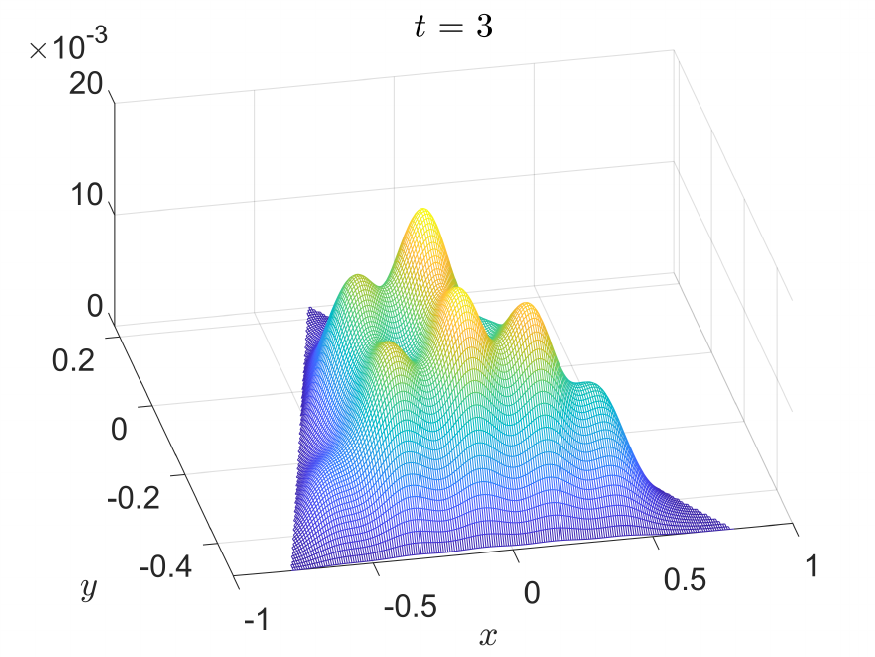}}\hspace*{-5pt}
  \subfigure[$t=6$]{\includegraphics[width=0.33\linewidth]{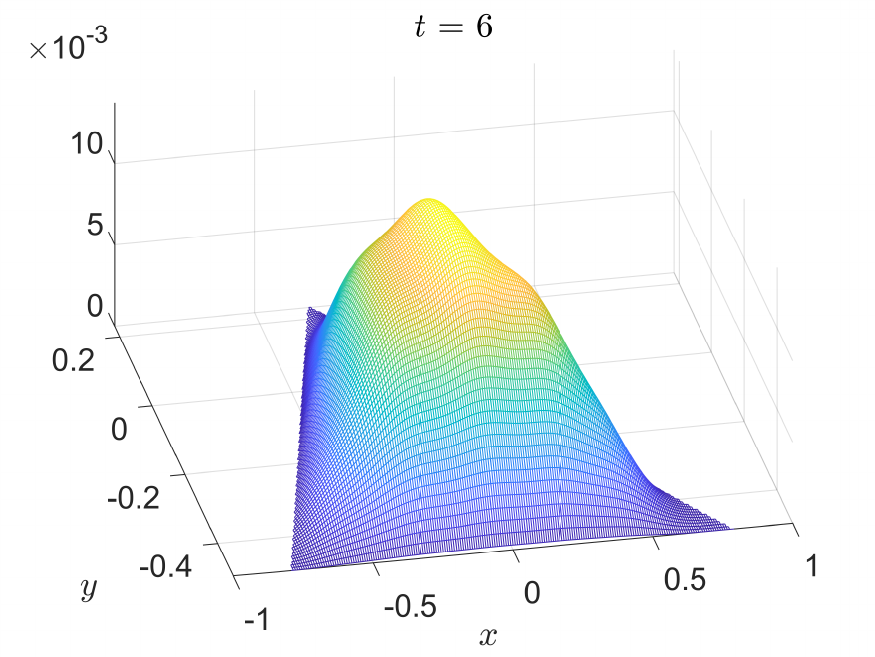}}
  \vspace*{-5pt}
  \caption{Viscous Burgers solution on the triangular domain at $t=0,3,6$, with
  $\nu=10^{-3}$, $\mu=1$, $\tau=10^{-2}$, and initial data~\eqref{eq:burg-init}.}
  \label{fig:burg-sol-tri}
\end{figure}

In \cref{fig:burg-sol-tri}, we depict the solutions of \eqref{eq:burg} defined on a triangular domain with coordinates $(-\tfrac45, -\tfrac12), (\tfrac45, -\tfrac12)$ and $(\tfrac45\cos\tfrac{5\pi}{8},\tfrac45\sin\tfrac{5\pi}{8}-\tfrac12)$.
Here we take $\nu=10^{-3}, \mu=1, \tau = 10^{-2}$ and $T=6.$
It can be observed that shock waves arising at the beginning fuse together.
Similar phenomena appear on other polygonal domains.

\section{Parameter selection and a one-dimensional prototype}\label{sec:opt-gam}
The pole parameters in the modified-Bessel solver are guided by the approximation
theory of lightning rational functions.  This theory is not, by itself, a convergence
proof for the Bessel ansatz~\eqref{eq:poi-app}. The numerical tests below examine
whether the same scaling laws are effective in the present setting.

For $f(x)=x^\alpha$ on $[0,1]$, $0<\alpha<1$, Stahl's result for the best type-$(N,N)$ rational approximation $r_N^*$ is
\begin{equation}\label{eq:rat-best}
  \lim_{N\to\infty}\re^{2\pi\sqrt{\alpha N}}
  \|f-r_N^*\|_{L^\infty(0,1)}
  =4^{1+\alpha}|\sin(\pi\alpha)|.
\end{equation}
Lightning approximations obtain root-exponential rates by prescribing poles that cluster near the branch point.  
Trefethen et al.~\cite{trefethen2021exponential} studied a tapered pole distribution of the form
\begin{equation}\label{eq:pole-trap}
  p_j=-\exp(-\sigma(\sqrt{N_1}-\sqrt j)), \quad 1\leq j\leq N_1.
\end{equation}
and established the root-exponential convergence rate $O(\exp(-2\pi\sqrt{\alpha N}))$ under the optimal parameter scaling
\begin{equation*}%
  \sigma_{\rm opt}=2\pi/\sqrt{\alpha}, \quad N_2=O(\sqrt{N_1}),
\end{equation*}
see also \cite{herremans2023resolution,xiang2025MC}.

For a sector with interior angle $\beta\pi$, $0<\beta<2$, and a singularity
$z^\alpha$, the corresponding tapered-pole prediction is
\begin{equation*}%
  \|z^\alpha-r_N\|_\infty = O (\exp(-\pi\sqrt{2(2-\beta)\alpha N}))
\end{equation*}
with the optimal parameter scaling $\sigma_{\rm opt} = \pi\sqrt{2(2-\beta)}/\sqrt{\alpha}$ (cf., e.g., \cite{xiang2024root,xiang2025root}).  
Notably, for the leading Laplace corner exponent $\alpha=1/\beta$, this gives
\begin{equation}\label{eq:optimal-sig}
  \sigma_{\rm opt}=\pi\sqrt{2(2-\beta)\beta}.
\end{equation}
Unless otherwise stated, the exterior pole clusters in our numerical
experiments use this corner-dependent parameter choice. Although these
formulas were originally developed for lightning rational approximations,
we employ them here as a parameter-selection strategy for \eqref{eq:poi-app},
with its effectiveness demonstrated by the numerical results below.

\begin{figure}[!t]
  \centering
  \subfigure[$\kappa^2=10$, $N_2\approx\tfrac12\sqrt{N_1}$]
  {\includegraphics[width=0.45\linewidth]{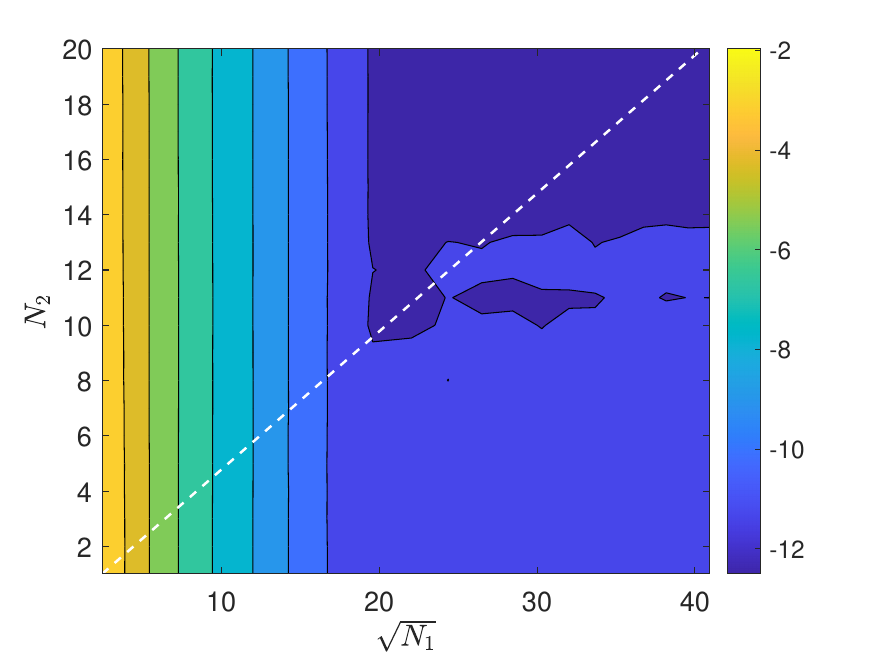}}
  \subfigure[$\kappa^2=1000$, $N_2\approx\tfrac45\sqrt{N_1}$]
  {\includegraphics[width=0.45\linewidth]{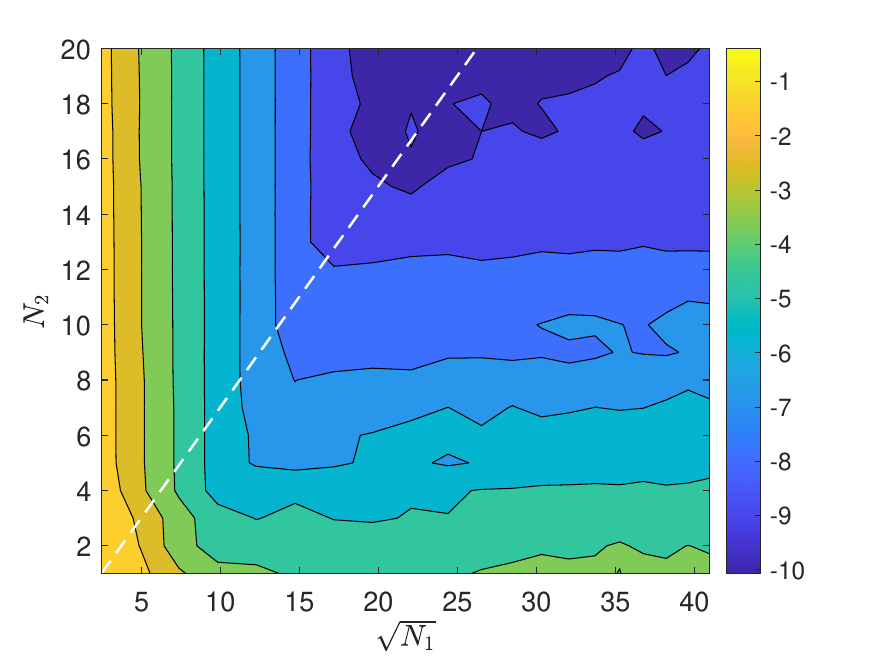}}
  \vspace*{-10pt}
  \caption{Boundary residuals for the L-shaped modified Helmholtz problem
  as $N_1$ and $N_2$ vary. The dashed white lines indicate the fitted
  relations $N_2=a\sqrt{N_1}+b$, corresponding to the fastest observed
  convergence with respect to the total number of coefficients
  $N=2N_1+2N_2+1$.}
  \label{fig:poitype-N1N2}
\end{figure}

We first investigate the influence of $N_1$ and $N_2$ for the L-shaped
modified Helmholtz problem in \cref{fig:poi-type-alp10}. The total number
of complex coefficients is $N=2N_1+2N_2+1.$
As shown in \cref{fig:poitype-N1N2}, the most efficient balance between
the two types of basis functions is well described by an empirical
square-root scaling $N_2\approx a\sqrt{N_1}+b,$
for both tested values of $\alpha$. This confirms the asymptotic relation
$N_2=O(\sqrt{N_1})$ suggested by the tapered pole distribution theory.

\begin{figure}[!t]
  \centering
  \includegraphics[width=0.45\linewidth]{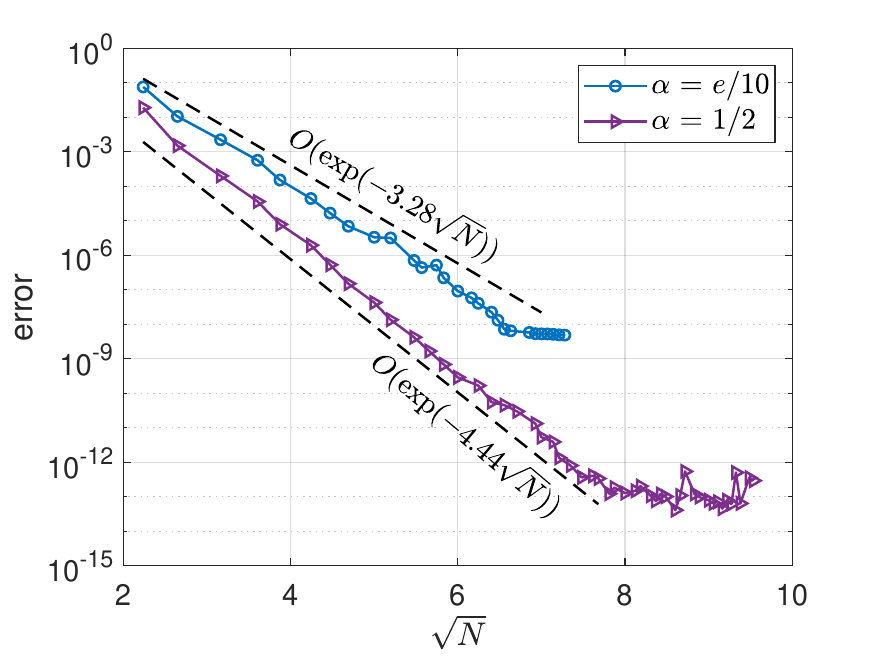}
  \includegraphics[width=0.45\linewidth]{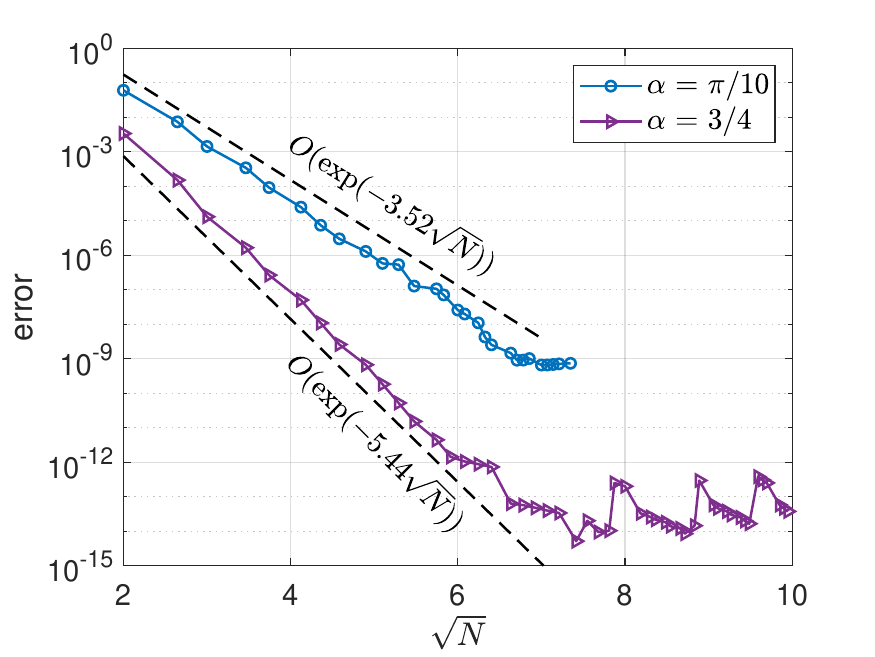}
  \vspace*{-10pt}
  \caption{Approximation errors of~\eqref{eq:app-bessel1} for $x^\alpha$ on
  $[0,1]$.  Dashed lines show the asymptotic slopes suggested by the lightning rational approximation theory: $\exp(-2\pi\sqrt{\alpha N})$.}
  \label{fig:lig-bessel0}
\end{figure}

As a simpler prototype, we consider the one-dimensional Bessel-type approximation
\begin{equation}\label{eq:app-bessel1}
  \widehat R_N(x) = \sum_{n=1}^{N_1}a_nK_1(x-x_n) + \sum_{m=0}^{N_2}b_mI_m(x-x_*),
  \qquad 0\leq x\leq1,
\end{equation}
where the exterior points $x_n<0$ are defined by \eqref{eq:pole-trap}, $x_*\in[0,1]$, and $N=N_1+N_2+1$.
In \cref{fig:lig-bessel0}, we compare the numerical errors for several
values of $\alpha$ with the reference slopes corresponding to the
predicted root-exponential convergence rate $\exp(-2\pi\sqrt{\alpha N}).$
The observed convergence behavior is consistent with the same
root-exponential scaling. A rigorous extension of the existing
rational-approximation convergence theory to the full two-dimensional
Bessel-type ansatz remains an interesting open problem.

\section{Extension to 3D modified Helmholtz equations}\label{sec:helm-3d}

Building on the key ideas for the two-dimensional modified Helmholtz equation developed in previous sections, we propose a three dimensional extension to deal with point singularities.
For the purposes, we first need some fundamental knowledge of spherical harmonic polynomials and refer to \cite{dai2013approximation} for more details.

\subsection{Spherical harmonic polynomials}
Harmonic polynomials are a set of polynomials that satisfy the Laplace equation, $\Delta Y = 0.$
Spherical harmonic functions form a subset of the harmonic polynomials, specifically those that are the restriction of homogeneous harmonic polynomials on the unit sphere.
Let $\mathcal H_l$ denote the space of spherical harmonic polynomials of degree $l$. The dimension of $\mathcal H_l$ is given by
\begin{equation*}
    d_l := {\rm dim} \mathcal H_l = \binom{l+2}{l} - \binom{l}{l-2} = 2l+1.
\end{equation*}
For any $l \in \mathbb N_0$, we denote by $\{Y_l^m: -l \le m \le l\}$ a set of spherical harmonic basis functions for $\mathcal H_l$. Then, the following property holds
\begin{align}\label{eq:prop-sphar}
    \Delta_0 Y = -l (l+1) Y, \quad \forall Y \in \mathcal H_l,\quad l=0,1,\ldots,
\end{align}
where $\Delta_0$ is the Laplace--Beltrami operator on the unit sphere.
Typically, the spherical harmonic polynomials can be explicitly expressed as \cite[(14.30.1)]{Olver2010Book}
\begin{align*}%
    Y_l^m(\xi) =  Y_l^m(\theta, \phi) 
    := (-1)^m \sqrt{ \frac{(2l+1) (l-m)!}{4\pi (l+m)!} } P_l^m(\cos\theta) \re^{\ri m \phi},
    \;\; -l \le m \le l,
\end{align*}
where $P_l^m(x)$ is the associated Legendre polynomials, defined by
\begin{equation*}
    P_l^m(x) := \frac{(-1)^m}{2^l l!} (1-x^2)^{m/2} \frac{\rd^{l+m}}{\rd x^{l+m}} (x^2-1)^l,
    \quad x \in [-1,1],
\end{equation*}
and $l$ and $m$ represent the degree and order of the associated Legendre polynomial, respectively, with $0 \le m \le l.$ 
Additionally, the associated Legendre polynomial satisfies the identity \cite[(14.9.3)]{Olver2010Book} %
\begin{equation*}
    P_l^{-m}(x) = (-1)^m \frac{(l-m)!}{(l+m)!} P_l^m(x),
\end{equation*}
which allows us to deduce the following symmetry for the spherical harmonic polynomials
\begin{equation*}
    Y_{l}^{-m}(\xi) = (-1)^m (Y_l^m)^*(\xi),
\end{equation*}
where the superscript $*$ denotes complex conjugation.

\subsection{Lightning solver for 3D modified Helmholtz equations}

Let us first consider the homogeneous equation defined on $\mathbb R^3$, i.e.,
\begin{align}\label{eq:helm-3d-gen2}
    -\Delta u + \kappa^2 u = 0,\quad \bs x:=(x,y,z) \in \mathbb R^3, \quad \kappa \in \mathbb R \setminus \{0\}.
\end{align}
In three-dimensional spherical polar coordinates, the Laplace operator $\Delta$ can be written as
\begin{equation*}
    \Delta = \frac{\rd^2}{\rd r^2} + \frac{2}{r} \frac{\rd}{\rd r} + \frac{1}{r^2} \Delta_0, \quad r = \|\bs x\|.
\end{equation*}

By representing the unknown $u$ as an expansion in spherical harmonic 
\begin{equation*}
    u(\bs x) = \sum_{l=0}^\infty \sum_{m=-l}^{l} \tilde u_l^m(r) Y_l^m(\xi), \quad \bs x = r\xi,\;\; r = \|\bs x\|, 
\end{equation*}
and recalling \eqref{eq:prop-sphar}, the modified Helmholtz equation \eqref{eq:helm-3d-gen2} becomes
\begin{align*}
    -\frac{1}{r^2} \partial_r (r^2 \partial_r \tilde u_l^m) + \frac{l(l+1)}{r^2} \tilde u_l^m + \kappa^2 \tilde u_l^m = 0,
\end{align*}
where we have used the property of spherical harmonics in \eqref{eq:prop-sphar}.
Next, we introduce a variable transformation $\eta = \kappa r$ and define $\tilde w_l^m(\eta) = r^{1/2} \tilde u_l^m(r)$. 
After substituting into the above formula, we consequently obtain the following equation
\begin{equation}\label{eq:helm-3d-gen3}
    \eta^2 \partial_\eta^2 \tilde w_l^m(\eta) + \eta \partial_\eta \tilde w_l^m(\eta) - (\eta^2 + \beta_l^2) \tilde w_l^m(\eta) = 0,
\end{equation}
where $\beta_l = l+\frac12.$

Notably, the second-order ordinary differential equation \eqref{eq:helm-3d-gen3} is exactly the Sturm--Liouville equation for the modified Bessel function.
Therefore, similar to \eqref{eq:poi-formal}, we can formally represent the solution to \eqref{eq:helm-3d-gen2} as
\begin{align}\label{eq:3d-sol10}\small
\begin{split}
    u(\bs x) = |\bs x - \bs x_*|^{-\frac12} \sum_{l=0}^\infty \sum_{m=-l}^{l} \big( a_l^m K_{\beta_l}(\kappa |\bs x - \bs x_*|) + b_l^m I_{\beta_l}(\kappa |\bs x - \bs x_*|) \big) Y_l^m(\xi)
\end{split}
\end{align}
for $\bs x \in \Omega$ and some $\bs x_* \in \mathbb R^3$, where $I_\nu$ and $K_\nu$ are the modified Bessel functions of the first and second kind, respectively.
According to the asymptotic behaviors of $K_\nu$ and $I_\nu$ from  \eqref{eq:asymp-IK}, we can modify \eqref{eq:3d-sol10} as
\begin{align}\label{eq:3d-sol2-mod}\small
\begin{split}
    u(\bs x) 
    &\approx \sum_{n=1}^{N_1} a_n |\bs x - \bs x_n|^{-\frac12} K_{\frac12}(\kappa |\bs x - \bs x_n|) Y_0^0(\xi_n) \\
    &\qquad + \sum_{l=0}^{N_2} \sum_{m=-l}^{l} b_l^m |\bs x - \bs x_*|^{-\frac12} I_{\beta_l} (\kappa |\bs x - \bs x_*|) Y_l^m(\xi_*), \quad \bs x \in \Omega
\end{split}
\end{align}
for $\{\bs x_n\}_{n=1}^{N_1} \not\subset \Omega $ and $\bs x_* \in \Omega.$
As noted earlier, we set $\bs x_n$ to be poles around the corners of the geometry $\Omega$, and $\bs x_*$ is chosen approximately at the center of $\Omega.$
In the last term of \eqref{eq:3d-sol2-mod}, we have, together with \eqref{eq:asymp-IK}, that 
\begin{equation*}\small
\begin{split}
|\bs x - \bs x_*|^{-\frac12} I_{\beta_l}(\kappa |\bs x - \bs x_*|) 
\to 
\frac{(\kappa/2)^{l+\frac12}}{\Gamma(l+3/2)} |\bs x - \bs x_*|^l,\quad \text{as $|\bs x - \bs x_*| \to 0$},
\end{split}
\end{equation*}
which is well-defined.
For simplicity, we focus on geometries with corners, such as a cone or the teardrop-shaped geometries.

\begin{remark}
For the inhomogeneous case, one can extend the source term $f$ and the particular solution $u_p$ to a ball with suitable boundary conditions. For numerical schemes, various efficient solvers can be employed $($see, for instance,  \cite{boulle2020computing} and references therein$)$.
\end{remark}

\subsection{Numerical experiments}
\begin{figure}[!h]
    \centering
    \subfigure[2D slice along the $y$-axis]{\includegraphics[width=0.42\linewidth]{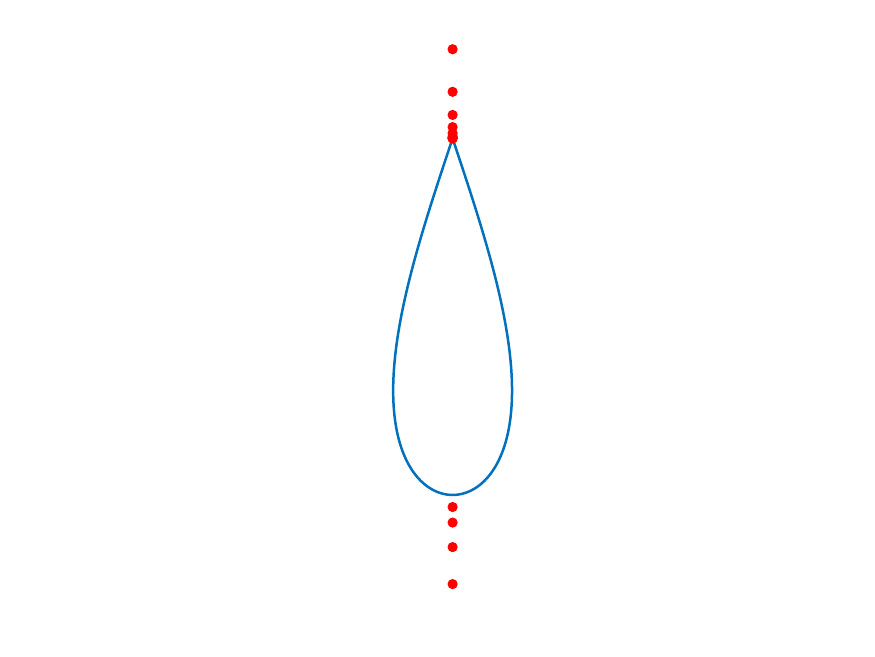}}
    \subfigure[3D visualization]{\includegraphics[width=0.45\linewidth]{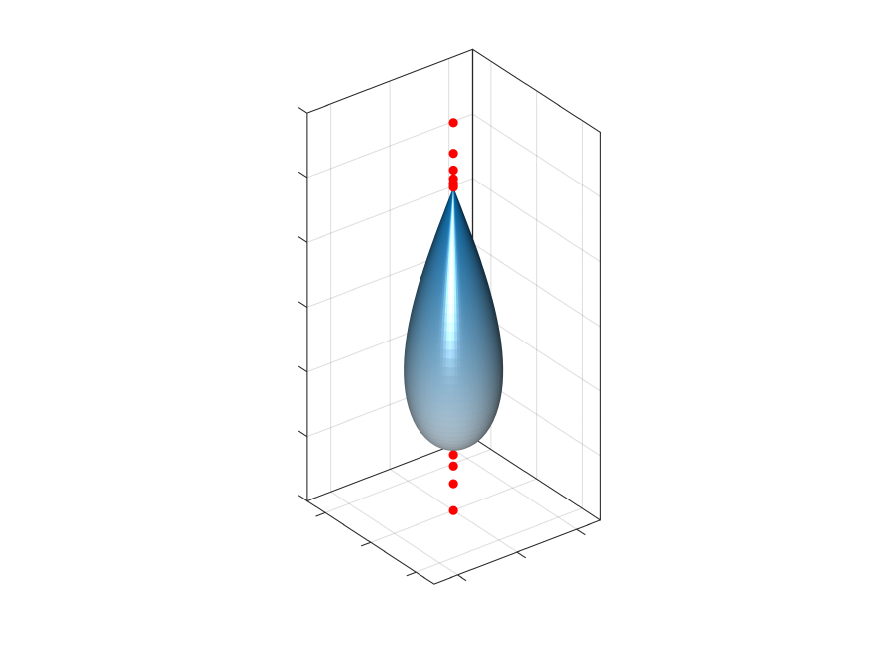}}
    \vspace*{-10pt}
    \caption{Illustration of the teardrop domain and the distribution of poles used in the lightning solver.}
    \label{fig:teardrop}
\end{figure}

To assess the efficiency of the proposed scheme \eqref{eq:3d-sol2-mod} numerically, we consider a teardrop-shaped geometry whose boundary is described by the parametric equation
\begin{equation*}
\bigl( \tfrac{\sqrt{3}}{3}\sin(2\pi\theta)\cos\varphi,\; \tfrac{\sqrt{3}}{3}\sin(2\pi\theta)\sin\varphi,\; -2\sin(\pi\theta) \bigr),
\end{equation*}
for $(\theta,\varphi) \in [0,1) \times [-\pi,\pi)$.
\cref{fig:teardrop} illustrates the geometry of the teardrop surface together with a two-dimensional slice along the $y$-axis.
The red dots indicate the locations of the poles, i.e., the points $\bs x_n$ in~\eqref{eq:3d-sol2-mod}.
In practice, we also place additional poles near the pseudo-singular point $(0,0,-2)$ along the $z$-axis in order to further enhance the efficiency and accuracy of the approximation.

For the numerical experiments, we select $M$ sampling points on the surface using a product grid.
Specifically, the grid is generated by taking
\texttt{theta = logspace(-16, log10(1/2), 300)}
(superposed with 50 Chebyshev points on $[0,1/2]$) in the $\theta$-direction, and 100 uniformly spaced points on $[-\pi,\pi)$ in the $\varphi$-direction.
As a result, the sampling points are strongly clustered near the singular point $(0,0,0)$.
To evaluate the numerical accuracy, the pointwise errors are computed on a second product grid consisting of medial points derived from the sampling grid.

The expansion coefficients $\{a_j, b_l^m\}$ are determined by solving the least-squares problem~\eqref{eq:poi-least}.
In the present setting, the entries of the matrix $\bs A$ take the form
\begin{equation*}
    \bs A_{j,n} =
    \begin{dcases}
        r_{j,n}^{-\frac12} K_{\frac12}(\kappa\, r_{j,n})\, Y_0^0(\xi),
        & 1 \le n \le N_1, \\[0.3em]
        r_{j,*}^{-\frac12} I_{\beta_l}(\kappa\, r_{j,*})\, Y_l^{-l+m}(\xi),
        & n = m + l^2 + N_1 + 1,
    \end{dcases}
\end{equation*}
for $1 \le j \le M$, $0 \le m \le 2l$, and $0 \le l \le N_2$.
Here, we choose $n_1$ poles with clustering parameter $\sigma = 4$ in~\eqref{eq:pole-trap} around the origin, and $n_2$ poles with $\sigma = 1$ around the point $(0,0,-2)$, so that $N_1 = n_1 + n_2$.

\begin{figure}
    \centering
    \includegraphics[width=0.8\linewidth]{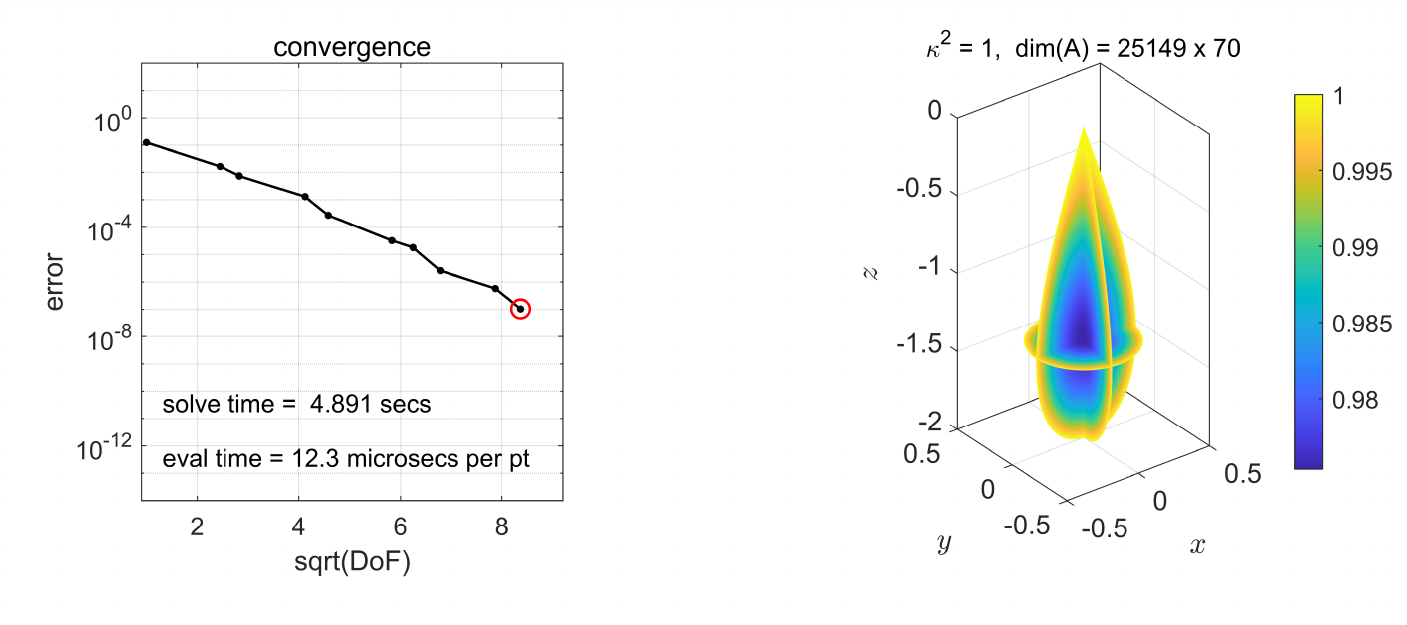}
    \includegraphics[width=0.8\linewidth]{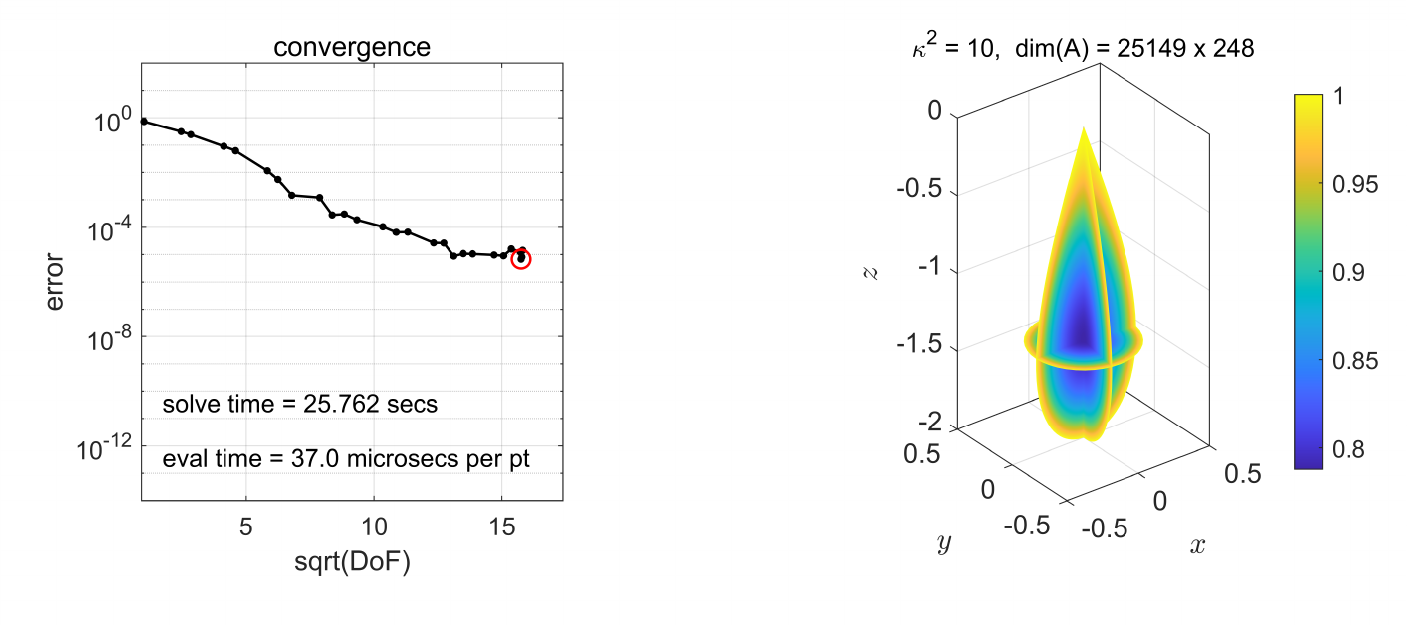}
    \vspace*{-10pt}
    \caption{Numerical solutions of the modified Helmholtz equation on a teardrop-shaped domain for $\kappa^2=1$ (top row) and $\kappa^2=10$ (bottom row).}
    \label{fig:lig-poitype-3d}
\end{figure}

The numerical results are presented in \cref{fig:lig-poitype-3d} with boundary conditions $h=1$.
The left panel shows the convergence behavior of the scheme~\eqref{eq:3d-sol2-mod}, while the right panel displays the computed solutions for $\kappa^2 = 1$ (top row) and $\kappa^2 = 10$ (bottom row), respectively.
We observe that, for $\kappa^2 = 1$, a solution with accuracy on the order of $10^{-6}$ is obtained in approximately 8 seconds.
For larger values of $\kappa$, achieving an accuracy of $10^{-5}$ requires additional computational time, as increasingly sharp boundary layers develop near the singular regions.
Similar behavior has also been observed in the two-dimensional case; see \cref{fig:poi-type-more}.

\section{Concluding remarks}\label{sec:conclud}
We developed a lightning architecture for bounded, inhomogeneous modified Helmholtz problems on two-dimensional cornered domains.  The homogeneous correction combines $K_1$ multipoles at tapered exterior pole clusters with a regular $I_m$ expansion, and its coefficients are determined from an oversampled, column-scaled least-squares problem.  
The inhomogeneous solver supplements this correction with a particular solution on an enclosing regular domain.  

The numerical tests show approximately root-exponential decay of independent boundary validation residuals on several polygonal domains and demonstrate the reuse of the solver in BDF discretizations of diffusion and viscous Burgers equations.  The one-dimensional Bessel prototype follows the slopes predicted by lightning rational approximation theory.  A complete convergence proof for the full two-dimensional Bessel ansatz, including the effects of boundary sampling and numerical rank truncation, remains an open problem.

We also described a three-dimensional extension to deal with isolated point singularities using a combination of  Bessel functions and spherical harmonic functions.
Although this approach can be in principle extended to deal with  three-dimensional  polyhedral edge and vertex singularities, how to efficiently solve the large scale least-squares problem required to resolve these singularities is particularly challenging.

\bibliographystyle{plain}
\bibliography{References}
\end{document}